# ANALYSIS OF PHASE TRANSITIONS IN THE MEAN-FIELD BLUME–EMERY–GRIFFITHS MODEL

By Richard S. Ellis[1], Peter T. Otto[1] and Hugo Touchette[2]

*University of Massachusetts, Union College and Queen Mary, University of London*

In this paper we give a complete analysis of the phase transitions in the mean-field Blume–Emery–Griffiths lattice-spin model with respect to the canonical ensemble, showing both a second-order, continuous phase transition and a first-order, discontinuous phase transition for appropriate values of the thermodynamic parameters that define the model. These phase transitions are analyzed both in terms of the empirical measure and the spin per site by studying bifurcation phenomena of the corresponding sets of canonical equilibrium macrostates, which are defined via large deviation principles. Analogous phase transitions with respect to the microcanonical ensemble are also studied via a combination of rigorous analysis and numerical calculations. Finally, probabilistic limit theorems for appropriately scaled values of the total spin are proved with respect to the canonical ensemble. These limit theorems include both central-limit-type theorems, when the thermodynamic parameters are not equal to critical values, and noncentral-limit-type theorems, when these parameters equal critical values.

**1. Introduction.** The Blume–Emery–Griffiths (BEG) model [4] is an important lattice-spin model in statistical mechanics. It is one of the few and certainly one of the simplest models known to exhibit, in its mean-field version, both a continuous, second-order phase transition and a discontinuous, first-order phase transition. Because of this property, the model has been studied extensively as a model of many diverse systems, including $He^3$-$He^4$ mixtures—the system for which Blume, Emery and Griffiths first devised

Received April 2004; revised January 2005.
[1]Supported by the NSF Grant DMS-02-02309.
[2]Supported by the Natural Sciences and Engineering Research Council of Canada and the Royal Society of London (Canada–UK Millennium Fellowship).

their model [4]—as well as solid-liquid-gas systems [18, 24, 25], microemulsions [23], semiconductor alloys [19] and electronic conduction models [17]. Phase diagrams for a class of models including the Blume–Emery–Griffiths model are discussed in [1], which lists additional work on this and related models. On a more theoretical level, the BEG model has also played an important role in the development of the renormalization-group theory of phase transitions of the Potts model; see [16, 20] for details and references.

As a model with a simple description but a relatively complicated phase transition structure, the BEG model continues to be of interest in modern statistical mechanical studies. In this paper we focus on the mean-field version of the BEG model or, equivalently, the BEG model on the complete graph on $n$ vertices. Our motivation for revisiting this model was initiated by a recent observation in [2, 3] that the BEG model on the complete graph has nonequivalent microcanonical and canonical ensembles, in the sense that it exhibits microcanonical equilibrium properties having no equivalent within the canonical ensemble. This observation is supported in [15] by numerical calculations both at the thermodynamic level, as in [2, 3], and at the level of equilibrium macrostates. In response to these earlier works, in this paper we address the phase transition behavior of the model by giving separate analyses of the structure of the sets of equilibrium macrostates for each of the two ensembles. Not only are our results consistent with the findings in [2, 3, 15], but also we rigorously prove for the first time a number of results that significantly generalize those found in these papers, where they were derived nonrigorously. For the canonical ensemble, full proofs of the structure of the set of equilibrium macrostates are provided. For the microcanonical ensemble, full proofs could not be attained. However, using numerical methods and following an analogous technique used in the canonical case, we also analyze the structure of the set of microcanonical equilibrium macrostates.

The BEG model that we consider is a spin-1 model defined on the complete graph on $n$ vertices $1, 2, \ldots, n$. The spin at site $j \in \{1, 2, \ldots, n\}$ is denoted by $\omega_j$, a quantity taking values in $\Lambda = \{-1, 0, 1\}$. The Hamiltonian for the BEG model is defined by

$$H_{n,K}(\omega) = \sum_{j=1}^{n} \omega_j^2 - \frac{K}{n}\left(\sum_{j=1}^{n} \omega_j\right)^2,$$

where $K > 0$ is a given parameter representing the interaction strength and $\omega = (\omega_1, \ldots, \omega_n) \in \Lambda^n$. The energy per particle is defined by

$$(1.1) \qquad h_{n,K}(\omega) = \frac{1}{n} H_{n,K}(\omega) = \frac{\sum_{j=1}^{n} \omega_j^2}{n} - K\left(\frac{\sum_{j=1}^{n} \omega_j}{n}\right)^2.$$

In order to analyze the phase transition behavior of the model, we first introduce the sets of equilibrium macrostates for the canonical ensemble and the microcanonical ensemble. As we will see, the canonical equilibrium macrostates solve a two-dimensional, unconstrained minimization problem, while the microcanonical equilibrium macrostates solve a dual, one-dimensional, constrained minimization problem. The definitions of these sets follow from large deviation principles derived for general models in [10]. In the particular case of the BEG model, they are consequences of the fact that the BEG-Hamiltonian can be written as a function of the empirical measures of the spin random variables and that, according to Sanov's theorem, the large deviation behavior of these empirical measures is governed by the relative entropy.

We use two innovations to analyze the structure of the set of canonical equilibrium macrostates. The first is to reduce to a one-dimensional problem the two-dimensional minimization problem that characterizes these macrostates. This is carried out by absorbing the noninteracting component of the energy per particle function into the prior measure, which is a product measure on configuration space. This manipulation allows us to express the canonical ensemble in terms of the empirical means, or spin per site $S_n/n = \sum_{j=1}^{n} \omega_j/n$, of the spin random variables. Doing so reduces the analysis of the BEG model to the analysis of a Curie–Weiss-type model [9] with single-site measures depending on $\beta$.

The analysis of the set of canonical equilibrium macrostates is further simplified by a second innovation. Because the thermodynamic parameter that defines the canonical ensemble is the inverse temperature $\beta$, a phase transition with respect to this ensemble is defined by fixing the Hamiltonian-parameter $K$ and varying $\beta$. Our analysis of the set of canonical equilibrium macrostates is based on a much more efficient approach that fixes $\beta$ and varies $K$. Proceeding in this way allows us to solve rigorously and in complete detail the reduced one-dimensional problem characterizing the equilibrium macrostates. We then extrapolate these results obtained by fixing $\beta$ and varying $K$ to physically relevant results that hold for fixed $K$ and varying $\beta$. These include a second-order, continuous phase transition and a first-order, discontinuous phase transition for different ranges of $K$.

For the microcanonical ensemble, we use a technique employed in [2] that absorbs the constraint into the minimizing function. This step allows us to reduce the constrained minimization problem defining the microcanonical equilibrium macrostates to another minimization problem that is more easily solved. Rigorous analysis of the reduced problem being limited, we rely mostly on numerical computations to complete our analysis of the set of equilibrium macrostates. Because the thermodynamic parameter defining the microcanonical ensemble is the energy per particle $u$, a phase transition with respect to this ensemble is defined by fixing $K$ and varying $u$. By analogy with the canonical case, our numerical analysis of the set of microcanonical equilibrium macrostates is based on a much more efficient approach that fixes $u$ and varies $K$. The analysis with respect to $K$ rather than $u$ allows us to solve in some detail the reduced problem characterizing the equilibrium macrostates. We then extrapolate these results obtained by fixing $u$ and varying $K$ to physically relevant results that hold for fixed $K$ and varying $u$. As in the case of the canonical ensemble, these include a second-order, continuous phase transition and a first-order, discontinuous phase transition for different ranges of $K$.

The contributions of this paper include a rigorous global analysis of the first- order phase transition in the canonical ensemble. Blume, Emery and Griffiths did a local analysis of the spin per site to show that their model exhibits a second-order phase transition for a range of values of $K$ and that, at a certain value of $K$, a tricritical point appears [4]; a similar study of a related model is carried out in [5, 6]. This tricritical point has the property that, for all smaller values of $K$, we are dealing with a first-order phase transition. Mathematically, the tricritical point marks the beginning of the failure of the local analysis; beyond this point, one has to resort to a global analysis of the spin per site. While the first-order phase transition has been studied numerically by several authors, the present paper gives the first rigorous global analysis.

Another contribution is that we analyze the phase transition for the canonical ensemble both in terms of the spin per site and the empirical measure. While all previous studies of the BEG model, except for [15], focused only on the spin per site, the analysis in terms of the empirical measure is the natural context for understanding equivalence and nonequivalence of ensembles [15].

A main consequence of our analysis is that the tricritical point—the critical value of the Hamiltonian parameter $K$ at which the model changes its phase transition behavior from second-order to first-order—differs in the two ensembles. Specifically, the tricritical point is smaller in the microcanonical ensemble than in the canonical ensemble. Therefore, there exists a range of values of $K$ such that the BEG model with respect to the canonical ensemble exhibits a first-order phase transition, while, with respect to the microcanonical ensemble, the model exhibits a second-order phase transition. As we discuss in Section 5, these results are consistent with the observation, seen numerically in [15], that there exists a subset of the microcanonical equilibrium macrostates that are not realized canonically. This observation implies that the two ensembles are nonequivalent at the level of equilibrium macrostates.

A final contribution of this paper is to present probabilistic limit theorems for appropriately scaled partial sums $S_n = \sum_{j=1}^{n} \omega_j$ with respect to the canonical ensemble. These limits follow from our work in Section 3 and known limit theorems for the Curie–Weiss model derived in [12, 14].

They include conditioned limit theorems when there are multiple equilibrium macrostates representing coexisting phases. In most cases the limits involve the central-limit-type scaling $n^{1/2}$ and convergence in distribution of $(S_n - n\tilde{z})/n^{1/2}$ to a normal random variable, where $\tilde{z}$ is an equilibrium macrostate. They also include the following two nonclassical cases, which hold for appropriate critical values of the parameters defining the canonical ensemble:

$$S_n/n^{3/4} \xrightarrow{\mathcal{D}} X \qquad \text{where } P\{X \in dx\} = \text{const} \cdot \exp[-\text{const} \cdot x^4] \, dx$$

and

$$S_n/n^{5/6} \xrightarrow{\mathcal{D}} X \qquad \text{where } P\{X \in dx\} = \text{const} \cdot \exp[-\text{const} \cdot x^6] \, dx.$$

As in the case of more complicated models, such as the Ising model, these nonclassical theorems signal the onset of a phase transition in the BEG model ([9], Section V.8). They are analogues of a result for the much simpler Curie–Weiss model ([9], Theorem V.9.5).

The outline of the paper is as follows. In Section 2, following the general procedure described in [10], we define the canonical ensemble, the microcanonical ensemble and the corresponding sets of equilibrium macrostates. In Section 3 the structure of the set of canonical equilibrium macrostates is studied. The initial analysis is carried out in Sections 3.2 and 3.3 at the level of the spin per site $S_n/n$ after the BEG model is written as a Curie–Weiss-type model in Section 3.1. In Sections 3.4 and 3.5 the information at the level of the spin per site is lifted to the level of the empirical measures of the spin random variables using the contraction principle, a main tool in the theory of large deviations. In Section 4 we present new theoretical insights into, and numerical results concerning, the structure of the set of microcanonical equilibrium macrostates. In Section 5 we discuss the implications of the results in the two previous sections concerning the nature of the phase transitions in the BEG model, which in turn is related to the phenomenon of ensemble nonequivalence at the level of equilibrium macrostates. Section 6 is devoted to probabilistic limit theorems for appropriately scaled sums $S_n$.

**2. Sets of equilibrium macrostates for the two ensembles.** The canonical and microcanonical ensembles are defined in terms of probability measures on a sequence of probability spaces $(\Lambda^n, \mathcal{F}_n)$. The configuration spaces $\Lambda^n$ consist of microstates $\omega = (\omega_1, \ldots, \omega_n)$ with each $\omega_j \in \Lambda = \{-1, 0, 1\}$, and $\mathcal{F}_n$ is the $\sigma$-field consisting of all subsets of $\Lambda^n$. We also introduce the $n$-fold product measure $P_n$ on $\Omega_n$ with identical one-dimensional marginals $\rho = \frac{1}{3}(\delta_{-1} + \delta_0 + \delta_1)$.

In terms of the energy per particle $h_{n,K}$ defined in (1.1), for each $n \in \mathbb{N}$, $\beta > 0$ and $K > 0$, the partition function is defined by

$$Z_n(\beta, K) = \int_{\Lambda^n} \exp[-n\beta h_{n,K}] \, dP_n.$$

For sets $B \in \mathcal{F}_n$, the canonical ensemble for the BEG model is the probability measure

$$P_{n,\beta,K}(B) = \frac{1}{Z_n(\beta,K)} \cdot \int_B \exp[-n\beta h_{n,K}]\, dP_n. \tag{2.1}$$

For $u \in \mathbb{R}$, $r > 0$, $K > 0$ and sets $B \in \mathcal{F}_n$, the microcanonical ensemble is the conditional probability measure

$$\begin{aligned} P_n^{u,r,K}(B) &= P_n\{B | h_{n,K} \in [u-r, u+r]\} \\ &= \frac{P_n\{B \cap \{h_{n,K} \in [u-r, u+r]\}\}}{P_n\{h_{n,K} \in [u-r, u+r]\}}. \end{aligned} \tag{2.2}$$

As we point out after (2.4), for appropriate values of $u$ and all sufficiently large $n$, the denominator is positive and, thus, $P_n^{u,r,K}$ is well defined.

The key to our analysis of the BEG model is to express both the canonical and the microcanonical ensembles in terms of the empirical measure $L_n$ defined for $\omega \in \Lambda^n$ by

$$L_n = L_n(\omega, \cdot) = \frac{1}{n} \sum_{j=1}^n \delta_{\omega_j}(\cdot).$$

$L_n$ takes values in $\mathcal{P}(\Lambda)$, the set of probability measures on $\Lambda = \{-1, 0, 1\}$. For $i \in \Lambda$, $L_n(\omega, \{i\})$ denotes the relative frequency of spins $\omega_j$ taking the value $i$. We rewrite $h_{n,K}$ as

$$\begin{aligned} h_{n,K}(\omega) &= \frac{\sum_{j=1}^n \omega_j^2}{n} - K\left(\frac{\sum_{j=1}^n \omega_j}{n}\right)^2 \\ &= \int_\Lambda y^2 L_n(\omega, dy) - K\left(\int_\Lambda y L_n(\omega, dy)\right)^2, \end{aligned}$$

and, for $\mu \in \mathcal{P}(\Lambda)$, we define

$$\begin{aligned} f_K(\mu) &= \int_\Lambda y^2 \mu(dy) - K\left(\int_\Lambda y \mu(dy)\right)^2 \\ &= (\mu_1 + \mu_{-1}) - K(\mu_1 - \mu_{-1})^2. \end{aligned} \tag{2.3}$$

The range of this function is the closed interval $[\min(1-K, 0), 1]$. In terms of $f_K$, we express $h_{n,K}$ in the form

$$h_{n,K}(\omega) = f_K(L_n(\omega)).$$

We appeal to the theory of large deviations to define the sets of canonical equilibrium macrostates and microcanonical equilibrium macrostates. Since any $\mu \in \mathcal{P}(\Lambda)$ has the form $\sum_{i=-1}^1 \mu_i \delta_i$, where $\mu_i \geq 0$ and $\sum_{i=-1}^1 \mu_i = 1$, $\mathcal{P}(\Lambda)$ can be identified with the set of probability vectors in $\mathbb{R}^3$. We topologize $\mathcal{P}(\Lambda)$ with the relative topology that this set inherits as a subset of $\mathbb{R}^3$. The relative entropy of $\mu \in \mathcal{P}(\Lambda)$ with respect to $\rho$ is defined by

$$R(\mu|\rho) = \sum_{i=-1}^{1} \mu_i \log(3\mu_i).$$

Sanov's theorem states that, with respect to the product measures $P_n$, the empirical measures $L_n$ satisfy the large deviation principle (LDP) on $\mathcal{P}(\Lambda)$ with rate function $R(\cdot|\rho)$ ([9], Theorem VIII.2.1). That is, for any closed subset $F$ of $\mathcal{P}(\Lambda)$, we have the large deviation upper bound

$$\limsup_{n\to\infty} \frac{1}{n} \log P_n\{L_n \in F\} \leq - \inf_{\mu \in F} R(\mu|\rho),$$

and, for any open subset $G$ of $\mathcal{P}(\Lambda)$, we have the large deviation lower bound

$$\limsup_{n\to\infty} \frac{1}{n} \log P_n\{L_n \in G\} \geq - \inf_{\mu \in G} R(\mu|\rho).$$

From the LDP for the $P_n$-distributions of $L_n$, we can derive the LDPs of $L_n$ with respect to the two ensembles $P_{n,\beta,K}$ and $P_n^{u,r,K}$. In order to state these LDPs, we introduce two basic thermodynamic functions, one associated with each ensemble. For $\beta > 0$ and $K > 0$, the basic thermodynamic function for the canonical ensemble is the canonical free energy

$$\varphi_K(\beta) = -\lim_{n\to\infty} \frac{1}{n} \log Z_n(\beta, K).$$

It follows from Theorem 2.4(a) in [10] that this limit exists for all $\beta > 0$ and $K > 0$ and is given by

$$\varphi_K(\beta) = \inf_{\mu \in \mathcal{P}(\Lambda)} \{R(\mu|\rho) + \beta f_K(\mu)\}.$$

For the microcanonical ensemble, the basic thermodynamic function is the microcanonical entropy

(2.4) $\qquad s_K(u) = -\inf\{R(\mu|\rho) : \mu \in \mathcal{P}(\Lambda), f_K(\mu) = u\}.$

Since $R(\mu|\rho) \geq 0$ for all $\mu$, $s_K(u) \in [-\infty, 0]$ for all $u$. We define $\mathrm{dom}\, s_K$ to be the set of $u \in \mathbb{R}$ for which $s_K(u) > -\infty$. Clearly, $\mathrm{dom}\, s_K$ coincides with the range of $f_K$ on $\mathcal{P}(\Lambda)$, which equals the closed interval $[\min(1-K, 0), 1]$. For $u \in \mathrm{dom}\, s_K$ and all sufficiently large $n$, the denominator in the second line of (2.2) is positive and, thus, the microcanonical ensemble $P_n^{u,r,K}$ is well defined ([10], Proposition 3.1).

The LDPs for $L_n$ with respect to the two ensembles are given in the next theorem. They are consequences of Theorems 2.4 and 3.2 in [10].

THEOREM 2.1. (a) *With respect to the canonical ensemble $P_{n,\beta,K}$, the empirical measures $L_n$ satisfy the LDP on $\mathcal{P}(\Lambda)$ with rate function*

$$I_{\beta,K}(\mu) = R(\mu|\rho) + \beta f_K(\mu) - \varphi_K(\beta). \tag{2.5}$$

(b) *With respect to the microcanonical ensemble $P_n^{u,r,K}$, the empirical measures $L_n$ satisfy the LDP on $\mathcal{P}(\Lambda)$, in the double limit $n \to \infty$ and $r \to 0$, with rate function*

$$I^{u,K}(\mu) = \begin{cases} R(\mu|\rho) + s_K(u), & \text{if } f_K(\mu) = u, \\ \infty, & \text{otherwise.} \end{cases} \tag{2.6}$$

For $\mu \in \mathcal{P}$ and $\varepsilon > 0$, we denote by $B(\mu, \varepsilon)$ the closed ball in $\mathcal{P}$ with center $\mu$ and radius $\varepsilon$. If $I_\beta(\mu) > 0$, then for all sufficiently small $\varepsilon > 0$, $\inf_{\nu \in B(\mu,\varepsilon)} I_\beta(\mu) > 0$. Hence, by the large deviation upper bound for $L_n$ with respect to the canonical ensemble, for all $\mu \in \mathcal{P}(\Lambda)$ satisfying $I_\beta(\mu) > 0$, all sufficiently small $\varepsilon > 0$ and all sufficiently large $n$,

$$P_{n,\beta,K}\{L_n \in B(\mu,\varepsilon)\} \leq \exp\left[-n\left(\inf_{\nu \in B(\mu,\varepsilon)} I_\beta(\nu)\right)\Big/2\right],$$

which converges to 0 exponentially fast. Consequently, the most probable macrostates $\nu$ solve $I_{\beta,K}(\nu) = 0$. It is therefore natural to define the set of canonical equilibrium macrostates to be

$$\begin{aligned} \mathcal{E}_{\beta,K} &= \{\nu \in \mathcal{P}(\Lambda) : I_{\beta,K}(\nu) = 0\} \\ &= \{\nu \in \mathcal{P}(\Lambda) : \nu \text{ minimizes } R(\nu|\rho) + \beta f_K(\nu)\}. \end{aligned} \tag{2.7}$$

Similarly, because of the large deviation upper bound for $L_n$ with respect to the microcanonical ensemble, it is natural to define the set of microcanonical equilibrium macrostates to be

$$\begin{aligned} \mathcal{E}^{u,K} &= \{\nu \in \mathcal{P}(\Lambda) : I^{u,K}(\nu) = 0\} \\ &= \{\nu \in \mathcal{P}(\Lambda) : \nu \text{ minimizes } R(\nu|\rho) \text{ subject to } f_K(\nu) = u\}. \end{aligned} \tag{2.8}$$

Each element $\nu$ in $\mathcal{E}_{\beta,K}$ and $\mathcal{E}^{u,K}$ has the form $\nu = \nu_{-1}\delta_{-1} + \nu_0\delta_0 + \nu_1\delta_1$ and describes an equilibrium configuration of the model in the corresponding ensemble. For $i = -1, 0, 1$, $\nu_i$ gives the asymptotic relative frequency of spins taking the value $i$.

In the next section we begin our study of the sets of equilibrium macrostates for the BEG model by analyzing $\mathcal{E}_{\beta,K}$.

**3. Structure of the set of canonical equilibrium macrostates.** In this section we give a complete description of the set $\mathcal{E}_{\beta,K}$ of canonical equilibrium macrostates for all values of $\beta$ and $K$. In contrast to all other studies of the model, which fix $K$ and vary $\beta$, we analyze the structure of $\mathcal{E}_{\beta,K}$ by fixing $\beta$ and varying $K$. As stated in Theorems 3.1 and 3.2, there exists a critical value of $\beta$, denoted by $\beta_c$ and equal to $\log 4$, such that $\mathcal{E}_{\beta,K}$ has two different forms for $0 < \beta \leq \beta_c$ and for $\beta > \beta_c$. Specifically, for fixed $0 < \beta \leq \beta_c$, $\mathcal{E}_{\beta,K}$ exhibits a continuous bifurcation as $K$ passes through a critical value $K_c^{(2)}(\beta)$, while for fixed $\beta > \beta_c$, $\mathcal{E}_{\beta,K}$ exhibits a discontinuous bifurcation as $K$ passes through a critical value $K_c^{(1)}(\beta)$. In Section 5 we show how to extrapolate this information to information concerning the phase transition behavior of the canonical ensemble for varying $\beta$: a continuous, second-order phase transition for all fixed, sufficiently large values of $K$ and a discontinuous, first-order phase transition for all fixed, sufficiently small values of $K$.

In terms of the uniform measure $\rho = \frac{1}{3}(\delta_{-1} + \delta_0 + \delta_1)$, we define

$$(3.1) \qquad \rho_\beta(d\omega_j) = \frac{1}{Z(\beta)} \cdot \exp(-\beta\omega_j^2)\rho(d\omega_j),$$

where $Z(\beta) = \int_\Lambda \exp(-\beta\omega_j^2)\rho(d\omega_j)$. The next two theorems give the form of $\mathcal{E}_{\beta,K}$ for $0 < \beta \leq \beta_c$ and for $\beta > \beta_c$. Theorem 3.1 will be proved in Section 3.5 as a consequence of Theorem 3.6, which is proved in Section 3.2.

THEOREM 3.1. *Define $\beta_c = \log 4$ and let $\rho_\beta$ be the measure defined in* (3.1). *For $0 < \beta \leq \beta_c$, the following conclusions hold:*

(a) *There exists a critical value $K_c^{(2)}(\beta) > 0$ defined in* (3.19) *and having the following properties:*

(i) *For $0 < K \leq K_c^{(2)}(\beta)$, $\mathcal{E}_{\beta,K} = \{\rho_\beta\}$.*

(ii) *For $K > K_c^{(2)}(\beta)$, there exist probability measures $\nu^+(\beta,K)$ and $\nu^-(\beta,K)$ in $\mathcal{P}(\Lambda)$ such that $\nu^+(\beta,K) \neq \nu^-(\beta,K) \neq \rho_\beta$ and $\mathcal{E}_{\beta,K} = \{\nu^+(\beta,K), \nu^-(\beta,K)\}$.*

(b) *If we write $\nu^+(\beta,K) = \nu_{-1}^+\delta_{-1} + \nu_0^+\delta_0 + \nu_1^+\delta_1$, then $\nu^-(\beta,K) = \nu_1^+\delta_{-1} + \nu_0^+\delta_0 + \nu_{-1}^+\delta_1$.*

(c) *For each choice of sign, $\nu^\pm(\beta,K)$ is a continuous function for $K > K_c^{(2)}(\beta)$, and as $K \to (K_c^{(2)}(\beta))^+$, $\nu^\pm(\beta,K) \to \rho_\beta$. Therefore, $\mathcal{E}_{\beta,K}$ exhibits a continuous bifurcation at $K_c^{(2)}(\beta)$.*

The continuous bifurcation described in part (c) of the theorem is an analogue of a second-order phase transition and explains the superscript 2 on the critical value $K_c^{(2)}(\beta)$. The next theorem shows that, for $\beta > \beta_c$, the set $\mathcal{E}_{\beta,K}$ exhibits a discontinuous bifurcation at a value $K_c^{(1)}(\beta)$. This analogue of a first-order phase transition explains the superscript 1 on the corresponding critical value $K_c^{(1)}(\beta)$. Theorem 3.2 will be proved in Section 3.5 as a consequence of Theorem 3.8, which is proved in Section 3.3. As we will see in the proof of the latter theorem, $K_c^{(1)}(\beta)$ is the unique zero of the function $A(K)$ defined in (3.31) for $K \geq K_1(\beta)$; $K_1(\beta)$ is specified in Lemma 3.9.

THEOREM 3.2. *Define $\beta_c = \log 4$ and let $\rho_\beta$ be the measure defined in* (3.1). *For $\beta > \beta_c$, the following conclusions hold:*

(a) *There exists a critical value $K_c^{(1)}(\beta) > 0$ having the following properties:*

   (i) *For $0 < K < K_c^{(1)}(\beta)$, $\mathcal{E}_{\beta,K} = \{\rho_\beta\}$.*

   (ii) *For $K = K_c^{(1)}(\beta)$, there exist probability measures $\nu^+(\beta, K_c^{(1)}(\beta))$ and $\nu^-(\beta, K_c^{(1)}(\beta))$ such that $\nu^+ \neq \nu^- \neq \rho_\beta$ and $\mathcal{E}_{\beta,K} = \{\rho_\beta, \nu^+(\beta, K_c^{(1)}(\beta)), \nu^-(\beta, K_c^{(1)}(\beta))\}$.*

   (iii) *For $K > K_c^{(1)}(\beta)$, there exist probability measures $\nu^+(\beta, K)$ and $\nu^-(\beta, K)$ such that $\nu^+(\beta, K) \neq \nu^-(\beta, K) \neq \rho_\beta$ and $\mathcal{E}_{\beta,K} = \{\nu^+(\beta, K), \nu^-(\beta, K)\}$.*

(b) *If we write $\nu^+(\beta, K) = \nu_{-1}^+ \delta_{-1} + \nu_0^+ \delta_0 + \nu_1^+ \delta_1$, then $\nu^-(\beta, K) = \nu_1^+ \delta_{-1} + \nu_0^+ \delta_0 + \nu_{-1}^+ \delta_1$.*

(c) *For each choice of sign, $\nu^\pm(\beta, K)$ is a continuous function for $K \geq K_c^{(1)}(\beta)$, and as $K \to (K_c^{(1)}(\beta))^+$, $\nu^\pm(\beta, K) \to \nu^\pm(\beta, K_c^{(1)}(\beta)) \neq \rho_\beta$. Therefore, $\mathcal{E}_{\beta,K}$ exhibits a discontinuous bifurcation at $K_c^{(1)}(\beta)$.*

We prove Theorems 3.1 and 3.2 in several steps. In the first step, carried out in Section 3.1, we absorb the noninteracting component of the energy per particle into the product measure of the canonical ensemble. This reduces the model to a Curie–Weiss-type model, which can be analyzed in terms of the empirical means $S_n/n = \sum_{j=1}^n \omega_j/n$. The structure of the set of canonical equilibrium macrostates for this Curie–Weiss-type model is analyzed in Section 3.2 for $0 < \beta \leq \beta_c$ and in Section 3.3 for $\beta > \beta_c$. In Section 3.4 we lift our results from the level of the empirical means up to the level of the empirical measures using the contraction principle, a main tool in the theory of large deviations. Finally, in Section 3.5 we derive Theorems 3.1 and 3.2 from the results derived in Section 3.2 for $0 < \beta \leq \beta_c$ and in Section 3.3 for $\beta > \beta_c$.

3.1. *Reduction to the Curie–Weiss model.* The first step in the proofs of Theorems 3.1 and 3.2 is to rewrite the canonical ensemble $P_{n,\beta,K}$ in the form of a Curie–Weiss-type model. We do this by absorbing the noninteracting component of the energy per particle $h_{n,K}$ into the product measure of $P_{n,\beta,K}$. Defining $S_n(\omega) = \sum_{j=1}^{n} \omega_j$, we write

$$\begin{aligned}
P_{n,\beta,K}(d\omega) &= \frac{1}{Z_n(\beta,K)} \cdot \exp[-n\beta h_{n,K}(\omega)] P_n(d\omega) \\
&= \frac{1}{Z_n(\beta,K)} \cdot \exp\left[-n\beta\left(\frac{\sum_{j=1}^{n}\omega_j^2}{n} - K\left(\frac{\sum_{j=1}^{n}\omega_j}{n}\right)^2\right)\right] P_n(d\omega) \\
&= \frac{1}{Z_n(\beta,K)} \cdot \exp\left[n\beta K\left(\frac{S_n(\omega)}{n}\right)^2\right] \prod_{j=1}^{n} \exp(-\beta\omega_j^2)\rho(d\omega_j) \\
&= \frac{(Z(\beta))^n}{Z_n(\beta,K)} \cdot \exp\left[n\beta K\left(\frac{S_n(\omega)}{n}\right)^2\right] P_{n,\beta}(d\omega).
\end{aligned}$$

In this formula $Z(\beta) = \int_\Lambda \exp(-\beta\omega_j^2)\rho(d\omega_j)$ and $P_{n,\beta}$ is the product measure on $\Lambda^n$ with identical one-dimensional marginals $\rho_\beta$ defined in (3.1).

We define

$$\tilde{Z}_n(\beta,K) = \int_{\Lambda^n} \exp\left[n\beta\left(\frac{S_n}{n}\right)^2\right] dP_{n,\beta}.$$

Since $P_{n,\beta,K}$ is a probability measure, it follows that

$$\tilde{Z}_n(\beta,K) = \frac{Z_n(\beta,K)}{[Z(\beta)]^n}$$

and, thus, that

(3.2) $$P_{n,\beta,K}(d\omega) = \frac{1}{\tilde{Z}_n(\beta,K)} \cdot \exp\left[n\beta K\left(\frac{S_n(\omega)}{n}\right)^2\right] P_{n,\beta}(d\omega).$$

By expressing the canonical ensemble in terms of the empirical means $S_n/n$, we have reduced the BEG model to a Curie–Weiss-type model. Cramér's theorem ([9], Theorem II.4.1) states that, with respect to the product measures $P_{n,\beta}$, $S_n/n$ satisfies the LDP on $[-1,1]$ with rate function

(3.3) $$J_\beta(z) = \sup_{t\in\mathbb{R}}\{tz - c_\beta(t)\}.$$

In this formula $c_\beta$ is the cumulant generating function defined by

(3.4) $$\begin{aligned}
c_\beta(t) &= \log \int_\Lambda \exp(t\omega_1)\rho_\beta(d\omega_1) \\
&= \log\left[\frac{1 + e^{-\beta}(e^t + e^{-t})}{1 + 2e^{-\beta}}\right].
\end{aligned}$$

$J_\beta$ is finite on the closed interval $[-1,1]$ and is differentiable on the open interval $(-1,1)$. This function is expressed in (3.3) as the Legendre–Fenchel transform of the finite, strictly convex, differentiable function $c_\beta$. By the theory of these transforms ([22], Theorem 25.1, [9], Theorem VI.5.3(d)), for each $z \in (-1, 1)$,

$$J'_\beta(z) = (c'_\beta)^{-1}(z). \tag{3.5}$$

From the LDP for $S_n/n$ with respect to $P_{n,\beta}$, Theorem 2.4 in [10] gives the LDP for $S_n/n$ with respect to the canonical ensemble written in the form (3.2).

THEOREM 3.3. *With respect to the canonical ensemble $P_{n,\beta,K}$ written in the form (3.2), the empirical means $S_n/n$ satisfy the LDP on $[-1, 1]$ with rate function*

$$\tilde{I}_{\beta,K} = J_\beta(z) - \beta K z^2 - \inf_{t \in \mathbb{R}}\{J_\beta(t) - \beta K t^2\}. \tag{3.6}$$

In Section 2 the canonical ensemble for the BEG model is expressed in terms of the empirical measures $L_n$. The corresponding set $\mathcal{E}_{\beta,K}$ of canonical equilibrium macrostates is defined as the set of probability measures $\nu \in \mathcal{P}(\Lambda)$ for which the rate function $I_{\beta,K}$ in the associated LDP satisfies $I_{\beta,K}(\nu) = 0$ [see (2.7)]. By contrast, in (3.2) the canonical ensemble is expressed in terms of the empirical means $S_n/n$. We now consider the set $\tilde{\mathcal{E}}_{\beta,K}$ of canonical equilibrium macrostates for the BEG model expressed in terms of the empirical means. Theorem 3.3 makes it natural to define $\tilde{\mathcal{E}}_{\beta,K}$ as the set of $z \in [-1, 1]$ for which the rate function in that theorem satisfies $\tilde{I}_{\beta,K}(z) = 0$. Since $z$ is a zero of this rate function if and only if $z$ minimizes $J_\beta(z) - \beta K z^2$, we have

$$\tilde{\mathcal{E}}_{\beta,K} = \{z \in [-1, 1] : z \text{ minimizes } J_\beta(z) - \beta K z^2\}. \tag{3.7}$$

As we will see in Theorem 3.13, each $z \in \tilde{\mathcal{E}}_{\beta,K}$ equals the mean of a corresponding measure $\nu \in \mathcal{E}_{\beta,K}$. Thus, each $z \in \tilde{\mathcal{E}}_{\beta,K}$ describes an equilibrium configuration of the model in terms of the specific magnetization, or the asymptotic average spin per site.

Although $J_\beta(z)$ can be computed explicitly, the expression is messy. Instead, we use an alternative characterization of $\tilde{\mathcal{E}}_{\beta,K}$ given in the next proposition to determine the points in that set. This proposition is a special case of Theorem A.1 in [7].

PROPOSITION 3.4. *For $z \in \mathbb{R}$, define*

$$G_{\beta,K}(z) = \beta K z^2 - c_\beta(2\beta K z). \tag{3.8}$$

*Then for each $\beta > 0$ and $K > 0$,*

$$\min_{|z| \le 1}\{J_\beta(z) - \beta K z^2\} = \min_{z \in \mathbb{R}}\{G_{\beta,K}(z)\}. \tag{3.9}$$

In addition, the global minimum points of $J_\beta(z) - \beta K z^2$ coincide with the global minimum points of $G_{\beta,K}$. As a consequence,

(3.10) $$\tilde{\mathcal{E}}_{\beta,K} = \{z \in \mathbb{R} : z \text{ minimizes } G_{\beta,K}(z)\}.$$

PROOF. The finite, convex function $f(z) = c_\beta(2\beta K z)/(2\beta K)$ has the Legendre–Fenchel transform

$$f^*(z) = \sup_{x \in \mathbb{R}}\{xz - f(x)\} = \begin{cases} J_\beta(z)/(2\beta K), & \text{for } |z| \leq 1, \\ \infty, & \text{for } |z| > 1. \end{cases}$$

We prove the proposition by showing the following three steps:

1. $\sup_{z \in \mathbb{R}}\{f(z) - z^2/2\} = \sup_{|z| \leq 1}\{z^2/2 - f^*(z)\}$.
2. Both suprema in step 1 are attained, the first for some $z \in \mathbb{R}$ and the second for some $z \in (-1, 1)$.
3. The global maximum points of $f(z) - z^2/2$ coincide with the global maximum points of $z^2/2 - f^*(z)$.

The proof uses three properties of Legendre–Fenchel transforms:

1. For all $z \in \mathbb{R}$, $f^{**}(z) = (f^*)^*(z)$ equals $f(z)$ ([9], Theorem VI.5.3(e)).
2. If for some $x \in \mathbb{R}$ and $z \in \mathbb{R}$, we have $z = f'(x)$, then $f(x) + f^*(z) = xz$ ([22], Theorem 25.1, [9], Theorem VI.5.3(c)). In particular, if $z = x$, then $f(x) + f^*(x) = x^2$.
3. If there exists $x \in (-1, 1)$ and $y \in \mathbb{R}$ such that

(3.11) $$f^*(z) \geq f^*(x) + y(z - x) \qquad \text{for all } z \in [-1, 1],$$

then $y = (f^*)'(x)$ ([22], Theorem 25.1). Hence, by properties 1 and 2,

$$f^*(x) + f^{**}(y) = f^*(x) + f(y) = xy.$$

In particular, if (3.11) is valid with $y = x$, then $f(x) + f^*(x) = x^2$.

Step 1 in the proof is a special case of Theorem C.1 in [8]. For completeness, we present the straightforward proof. Let $M = \sup_{z \in \mathbb{R}}\{f(z) - z^2/2\}$. Since for any $|z| \leq 1$ and $x \in \mathbb{R}$

$$f^*(z) + M \geq xz - f(x) + M \geq xz - x^2/2,$$

we have

$$f^*(z) + M \geq \sup_{x \in \mathbb{R}}\{xz - x^2/2\} = z^2/2.$$

It follows that $M \geq z^2/2 - f^*(z)$ and thus that $M \geq \sup_{|z| \leq 1}\{z^2/2 - f^*(z)\}$. To prove the reverse inequality, let $N = \sup_{|z| \leq 1}\{z^2/2 - f^*(z)\}$. Then for any $z \in \mathbb{R}$ and $|x| \leq 1$,

$$z^2/2 + N \geq xz - x^2/2 + N \geq xz - f^*(x).$$

Since $f^*(x) = \infty$ for $|x| > 1$, it follows from property 1 that
$$z^2/2 + N \geq \sup_{|x| \leq 1} \{xz - f^*(x)\} = f(z)$$
and thus that $N \geq \sup_{z \in \mathbb{R}} \{f(z) - z^2/2\}$. This completes the proof of step 1.

Since $f(z) \sim |z|$ as $z \to \infty$, $f(z) - z^2/2$ attains its supremum over $\mathbb{R}$. Since $z^2/2 - f^*(z)$ is continuous and $\lim_{|z| \to 1} (f^*)'(z) = \infty$, $z^2/2 - f^*(z)$ attains its supremum over $[-1, 1]$ in the open interval $(-1, 1)$. This completes the proof of step 2.

We now prove that the global maximum points of the two functions coincide. Let $x$ be any point in $\mathbb{R}$ at which $f(z) - z^2/2$ attains its supremum. Then $x = f'(x)$, and so by the second assertion in property 2, $f(x) + f^*(x) = x^2$. The point $x$ lies in $(-1, 1)$ because the range of $f'(z) = c'_\beta(2\beta K z)$ equals $(-1, 1)$. Step 1 now implies that
$$\sup_{|z| \leq 1} \{z^2/2 - f^*(z)\} = \sup_{z \in \mathbb{R}} \{f(z) - z^2/2\}$$
$$= f(x) - x^2/2 = x^2/2 - f^*(x).$$
We conclude that $z^2/2 - f^*(z)$ attains its supremum at $x \in (-1, 1)$.

Conversely, let $x$ be any point in $(-1, 1)$ at which $z^2/2 - f^*(z)$ attains its supremum. Then for any $z \in [-1, 1]$,
$$x^2/2 - f^*(x) \geq z^2/2 - f^*(z).$$
It follows that, for any $z \in [-1, 1]$,
$$f^*(z) \geq f^*(x) + (z^2 - x^2)/2 \geq f^*(x) + x(z - x).$$
The second assertion in property 3 implies that $f^*(x) + f(x) = x^2$, and, in conjunction with step 1, this in turn implies that
$$\sup_{z \in \mathbb{R}} \{f(z) - z^2/2\} = \sup_{|z| \leq 1} \{z^2/2 - f^*(z)\}$$
$$= x^2/2 - f^*(x) = f(x) - x^2/2.$$
We conclude that $f(z) - z^2/2$ attains its supremum at $x$. This completes the proof of the proposition. $\square$

Proposition 3.4 states that $\tilde{\mathcal{E}}_{\beta,K}$ consists of the global minimum points of $G_{\beta,K}(z) = \beta K z^2 - c_\beta(2\beta K z)$. In order to simplify the minimization problem, we make the change of variables $z \to z/(2\beta K)$ in $G_{\beta,K}$, obtaining the new function

(3.12) $$F_{\beta,K}(z) = G_{\beta,K}\left(\frac{z}{2\beta K}\right) = \frac{z^2}{4\beta K} - c_\beta(z).$$

Proposition 3.4 gives the alternative characterization of $\tilde{\mathcal{E}}_{\beta,K}$ to be

$$\tilde{\mathcal{E}}_{\beta,K} = \left\{ \frac{w}{2\beta K} \in \mathbb{R} : w \text{ minimizes } F_{\beta,K}(w) \right\}. \tag{3.13}$$

We use $F_{\beta,K}$ to analyze $\tilde{\mathcal{E}}_{\beta,K}$ because the second term of $F_{\beta,K}$ contains only the parameter $\beta$, while both terms in $G_{\beta,K}$ contain both parameters $\beta$ and $K$. In order to analyze the structure of $\tilde{\mathcal{E}}_{\beta,K}$, we take advantage of the simpler form of $F_{\beta,K}$ by fixing $\beta$ and varying $K$. This innovation makes the analysis of $\tilde{\mathcal{E}}_{\beta,K}$ much more efficient than in previous studies. Our goal is prove that the elements of $\tilde{\mathcal{E}}_{\beta,K}$ change continuously with $K$ for all $0 < \beta \leq \beta_c = \log 4$ (Theorem 3.1) and have a discontinuity at $K_c^{(1)}$ for all $\beta > \beta_c$ (Theorem 3.2).

In order to determine the minimum points of $F_{\beta,K}$ and, thus, the points in $\tilde{\mathcal{E}}_{\beta,K}$, we study the derivative

$$F'_{\beta,K}(w) = \frac{w}{2\beta K} - c'_\beta(w). \tag{3.14}$$

$F'_{\beta,K}(w)$ consists of a linear part $w/(2\beta K)$ and a nonlinear part $c'_\beta(w)$. As we will see in Sections 3.2 and 3.3, the basic mechanism underlying the change in the bifurcation behavior of $\tilde{\mathcal{E}}_{\beta,K}$ is the change in the concavity behavior of $c'_\beta(w)$ for $0 < \beta \leq \beta_c$ versus $\beta > \beta_c$, which is the subject of the next theorem. A related phenomenon was observed in [11], Theorem 1.2(b), and in [13], Theorem 4, in the context of work on the Griffiths–Hurst–Sherman correlation inequality for models of ferromagnets; this inequality is used to show the concavity of the specific magnetization as a function of the external field.

THEOREM 3.5. *For $\beta > \beta_c = \log 4$, define*

$$w_c(\beta) = \cosh^{-1}(\tfrac{1}{2}e^\beta - 4e^{-\beta}) \geq 0. \tag{3.15}$$

*The following conclusions hold:*

(a) *For $0 < \beta \leq \beta_c$, $c'_\beta(w)$ is strictly concave for $w > 0$.*

(b) *For $\beta > \beta_c$, $c'_\beta(w)$ is strictly convex for $0 < w < w_c(\beta)$ and $c'_\beta(w)$ is strictly concave for $w > w_c(\beta)$.*

PROOF. (a) We show that for all $0 < \beta \leq \beta_c$, $c'''_\beta(w) < 0$ for all $w > 0$. A short calculation yields

$$c'''_\beta(w) = \frac{[2e^{-\beta} \sinh w][1 - 2e^{-\beta} \cosh w - 8e^{-2\beta}]}{[1 + 2e^{-\beta} \cosh w]^3}. \tag{3.16}$$

Since $2e^{-\beta} \sinh w$ and $1 + 2e^{-\beta} \cosh w$ are positive for $w > 0$, $c'''_\beta(w) < 0$ for $w > 0$ if and only if

$$1 - 2e^{-\beta} \cosh w - 8e^{-2\beta} < 0 \text{ for } w > 0.$$

The inequality $\cosh w > 1$ for $w > 0$ implies that

$$[1 - 2e^{-\beta}\cosh w - 8e^{-2\beta}] < [1 - 2e^{-\beta} - 8e^{-2\beta}]$$
$$= (1 - 4e^{-\beta})(1 + 2e^{-\beta}) \quad \text{for all } w > 0.$$

Therefore, for all $0 < \beta \leq \log 4$, $c'''_\beta(w) < 0$ for $w > 0$.

(b) Fixing $\beta > \beta_c$, we determine the critical value $w_c(\beta)$ such that $c'_\beta(w)$ is strictly convex for $0 < w < w_c(\beta)$ and strictly concave for $w > w_c(\beta)$. From the expression for $c'''_\beta(w)$ in (3.16), $c'''_\beta(w) > 0$ for $w > 0$ if and only if $(1 - 2e^{-\beta}\cosh w - 8e^{-2\beta}) > 0$ for $w > 0$. Therefore, $c'_\beta(w)$ is strictly convex for

$$0 < w < \cosh^{-1}(\tfrac{1}{2}e^\beta - 4e^{-\beta}).$$

On the other hand, since $c'''_\beta(w) < 0$ for $w > 0$ if and only if $(1 - 2e^{-\beta}\cosh w - 8e^{-2\beta}) < 0$ for $w > 0$, we conclude that $c'_\beta(w)$ is strictly concave for

$$w > \cosh^{-1}(\tfrac{1}{2}e^\beta - 4e^{-\beta}).$$

This completes the proof of part (b). □

The concavity description of $c'_\beta$ stated in Theorem 3.5 allows us to find the global minimum points of $F_{\beta,K}$ and thus the points in $\tilde{\mathcal{E}}_{\beta,K}$ for all values of the parameters $\beta$ and $K$. We carry this out in the next two sections, first for $0 < \beta \leq \beta_c$ and then for $\beta > \beta_c$. In Section 3.4 we use this information to give the structure of the set $\mathcal{E}_{\beta,K}$ of canonical equilibrium macrostates defined in (2.7).

3.2. *Description of $\tilde{\mathcal{E}}_{\beta,K}$ for $0 < \beta \leq \beta_c$.* In Theorem 3.1 we state the structure of the set $\mathcal{E}_{\beta,K}$ of canonical equilibrium macrostates for the BEG model with respect to the empirical measures when $0 < \beta \leq \beta_c = \log 4$. The main theorem in this section, Theorem 3.6, does the same for the set $\tilde{\mathcal{E}}_{\beta,K}$, which has been shown to have the alternative characterization

(3.17) $$\tilde{\mathcal{E}}_{\beta,K} = \left\{ \frac{w}{2\beta K} \in \mathbb{R} : w \text{ minimizes } F_{\beta,K}(w) \right\}.$$

We recall that $F_{\beta,K}(w) = w^2/(4\beta K) - c_\beta(w)$, where $c_\beta$ is defined in (3.4). In Section 3.4 we will prove that there exists a one-to-one correspondence between $\tilde{\mathcal{E}}_{\beta,K}$ and $\mathcal{E}_{\beta,K}$. In Section 3.5 we will use this fact to fully describe the latter set for all $0 < \beta \leq \beta_c$ and $K > 0$.

According to part (a) of Theorem 3.5, for $0 < \beta \leq \beta_c$, $c'_\beta(w)$ is strictly concave for $w > 0$. As a result, the study of $\tilde{\mathcal{E}}_{\beta,K}$ is similar to the study of the equilibrium macrostates for the classical Curie–Weiss model as given in Section IV.4 of [9]. Following the discussion in that section, we first use a graphical argument to motivate the continuous bifurcation exhibited by $\tilde{\mathcal{E}}_{\beta,K}$ for $0 < \beta \leq \beta_c$. A detailed statement is given in Theorem 3.6.

Minimum points of $F_{\beta,K}$ satisfy $F'_{\beta,K}(w) = 0$, which can be rewritten as

$$(3.18) \qquad \frac{w}{2\beta K} = c'_\beta(w).$$

Since the slope of the function $w \mapsto w/(2\beta K)$ is $1/(2\beta K)$, the nature of the solutions of (3.18) depends on whether

$$c''_\beta(0) \leq \frac{1}{2\beta K} \quad \text{or} \quad 0 < \frac{1}{2\beta K} < c''_\beta(0).$$

This motivates the definition of the critical value

$$(3.19) \qquad K_c^{(2)}(\beta) = \frac{1}{2\beta c''_\beta(0)} = \frac{1}{4\beta e^{-\beta}} + \frac{1}{2\beta}.$$

We use the same notation here as for the critical value in Theorem 3.1 because, as we will later prove, the continuous bifurcation in $K$ exhibited by both sets $\mathcal{E}_{\beta,K}$ and $\tilde{\mathcal{E}}_{\beta,K}$ occur at the same value $K_c^{(2)}(\beta)$.

We illustrate the minimum points of $F_{\beta,K}$ graphically in Figure 1 for $\beta = 1$. For three ranges of values of $K$, this figure depicts the two components of $F'_{\beta,K}$: the linear component $w/(2\beta K)$ and the nonlinear component $c'_\beta(w)$. Figure 1(a) corresponds to $0 < K < K_c^{(2)}(\beta)$. Since $c''_\beta(0) = 1/(2\beta K_c^{(2)}(\beta))$, for $0 < K < K_c^{(2)}(\beta)$, the two components of $F'_{\beta,K}$ intersect at only the origin, and, thus, $F_{\beta,K}$ has a unique global minimum point at $w = 0$. Figure 1(b) corresponds to $K = K_c^{(2)}(\beta)$. In this case the two components of $F'_{\beta,K}$ are tangent at the origin, and again $F_{\beta,K}$ has a unique global minimum point at $w = 0$. Figure 1(c) corresponds to $K > K_c^{(2)}(\beta)$. For such $K$, the global minimum points of $F_{\beta,K}$ are symmetric nonzero points $w = \pm\tilde{w}(\beta,K)$, $\tilde{w}(\beta,K) > 0$.

Figures 1(a) and 1(c) give similar information as Figures IV.3(b) and IV.3(d) in [9], which depict the phase transition in the Curie–Weiss model. In these two sets of figures the functions being graphed are Legendre–Fenchel transforms of each other.

The graphical information just obtained concerning the global minimum points of $F_{\beta,K}$ for $0 < \beta \leq \beta_c$ motivates the form of $\tilde{\mathcal{E}}_{\beta,K}$ stated in the next theorem. The positive quantity $\tilde{z}(\beta,K)$ equals $\tilde{w}(\beta,K)/(2\beta K)$; $\tilde{w}(\beta,K)$ is the unique positive global minimum point of $F_{\beta,K}$ for $K > K_c^{(2)}(\beta)$, the existence of which is proved in Lemma 3.7. According to part (c) of the theorem, $\tilde{z}(\beta,K)$ is a continuous function for $K > K_c^{(2)}(\beta)$, and as $K \to (K_c^{(2)}(\beta))^+$, $\tilde{z}(\beta,K)$ converges to 0. As a result, the bifurcation exhibited by $\tilde{\mathcal{E}}_{\beta,K}$ at $K_c^{(2)}(\beta)$ is continuous.

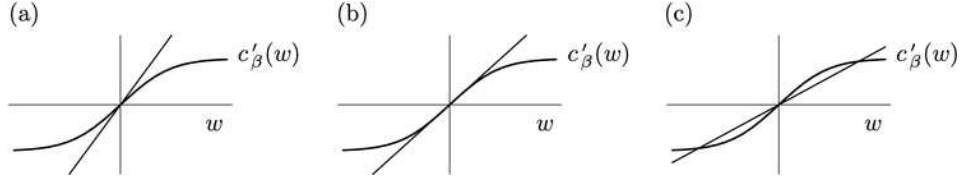

FIG. 1. *Continuous bifurcation for $\beta = 1$.* (a) $K < K_c^{(2)}(\beta)$, (b) $K = K_c^{(2)}(\beta)$, (c) $K > K_c^{(2)}(\beta)$.

THEOREM 3.6. *Define $\tilde{\mathcal{E}}_{\beta,K}$ by* (3.7); *equivalently,*

$$\tilde{\mathcal{E}}_{\beta,K} = \left\{ \frac{w}{2\beta K} \in \mathbb{R} : w \text{ minimizes } F_{\beta,K}(w) \right\}.$$

*For all $0 < \beta \leq \beta_c$, the critical value $K_c^{(2)}(\beta) = 1/(2\beta c_\beta''(0))$ has the following properties:*

(a) *For $0 < K \leq K_c^{(2)}(\beta)$, $\tilde{\mathcal{E}}_{\beta,K} = \{0\}$.*

(b) *For $K > K_c^{(2)}(\beta)$, there exists a positive number $\tilde{z}(\beta, K)$ such that $\tilde{\mathcal{E}}_{\beta,K} = \{\pm \tilde{z}(\beta, K)\}$.*

(c) *$\tilde{z}(\beta, K)$ is a strictly increasing continuous function for $K > K_c^{(2)}(\beta)$, and as $K \to (K_c^{(2)}(\beta))^+$, $\tilde{z}(\beta, K) \to 0$. Therefore, $\tilde{\mathcal{E}}_{\beta,K}$ exhibits a continuous bifurcation at $K_c^{(2)}(\beta)$.*

The proof of the theorem depends on the next lemma, in which we show that, for $K > K_c^{(2)}(\beta)$, $F_{\beta,K}$ has a unique positive global minimum at a point $\tilde{w}(\beta, K)$.

LEMMA 3.7. *For $0 < \beta \leq \beta_c = \log 4$, define $F_{\beta,K}$ by* (3.12). *The following conclusions hold:*

(a) *For each $K > K_c^{(2)}(\beta)$, $F_{\beta,K}$ has a critical point $\tilde{w}(\beta, K) > 0$ satisfying*

$$F'_{\beta,K}(\tilde{w}(\beta, K)) = 0 \quad \text{and} \quad F''_{\beta,K}(\tilde{w}(\beta, K)) > 0.$$

(b) *For each $K > K_c^{(2)}(\beta)$, $F_{\beta,K}$ has unique nonzero global minimum points at $w = \pm \tilde{w}(\beta, K)$.*

(c) *The points $\{\tilde{w}(\beta, K), K > K_c^{(2)}(\beta)\}$ span the positive real line; that is, for each $x > 0$, there exists $K > K_c^{(2)}(\beta)$ such that $x = \tilde{w}(\beta, K)$.*

PROOF. (a) For any $K > K_c^{(2)}(\beta)$, we have

$$(3.20) \qquad \frac{1}{2\beta K} < \frac{1}{2\beta K_c^{(2)}(\beta)} = c_\beta''(0) = \lim_{w \to 0} \frac{c_\beta'(w)}{w}.$$

Since $c'_\beta$ is continuous, for sufficiently small $w > 0$, we have $w/(2\beta K) < c'_\beta(w)$ and, thus, $F'_{\beta,K}(w) < 0$. On the other hand, $|c'_\beta(w)| < 1$ for all $w$ and, therefore, $\lim_{w\to\infty} F'_{\beta,K}(w) = \infty$. It follows that $F'_{\beta,K}(w) > 0$ for sufficiently large $w > 0$. Consequently, by the continuity of $F'_{\beta,K}$, there exists at least one positive critical point of $F_{\beta,K}$; the analyticity of $F_{\beta,K}$ implies that $F_{\beta,K}$ has at most finitely many critical points. Denote by $\tilde{w}(\beta, K) > 0$ the smallest positive critical point of $F_{\beta,K}$.

We now prove that $F''_{\beta,K}(\tilde{w}(\beta, K)) > 0$. Since $F'_{\beta,K}(\tilde{w}(\beta, K)) = 0$, the mean value theorem yields the existence of $\alpha \in (0, \tilde{w}(\beta, K))$ such that
$$c''_\beta(\alpha) = \frac{c'_\beta(\tilde{w}(\beta, K))}{\tilde{w}(\beta, K)} = \frac{1}{2\beta K}.$$

By part (a) of Theorem 3.5, since $\alpha < \tilde{w}(\beta, K)$, it follows that $c''_\beta(\alpha) > c''_\beta(\tilde{w}(\beta, K))$ and thus that

$$(3.21) \quad F''_{\beta,K}(\tilde{w}(\beta, K)) = \frac{1}{2\beta K} - c''_\beta(\tilde{w}(\beta, K)) > \frac{1}{2\beta K} - c''_\beta(\alpha) = 0.$$

This completes the proof of part (a).

(b) For any $w > \tilde{w}(\beta, K)$, the strict concavity of $c'_\beta(w)$ for $w > 0$ [Theorem 3.5(a)] implies that $c''_\beta(w) < c''_\beta(\tilde{w}(\beta, K))$. Therefore, by (3.21), we have
$$F''_{\beta,K}(w) = \frac{1}{2\beta K} - c''_\beta(w)$$
$$> \frac{1}{2\beta K} - c''_\beta(\tilde{w}(\beta, K)) = F''_{\beta,K}(\tilde{w}(\beta, K)) > 0.$$

Thus, $F'_{\beta,K}$ is strictly increasing for $w > \tilde{w}(\beta, K)$. This property allows us to conclude that $\tilde{w}(\beta, K)$ is the unique positive critical point and the unique positive local minimum point of $F_{\beta,K}$. By symmetry, $F_{\beta,K}$ has a unique negative local minimum point at $w = -\tilde{w}(\beta, K)$. In addition, as shown in (3.20), for any $K > K_c^2(\beta)$, we have $F''_{\beta,K}(0) = 1/(2\beta K) - c''_\beta(0) < 0$. Since $\lim_{|w|\to\infty} F_{\beta,K}(w) = \infty$, we conclude that $\pm\tilde{w}(\beta, K)$ are the unique global minimum points of $F_{\beta,K}$.

(c) Given $x > 0$, define the positive number $K_x = x/(2\beta c'_\beta(x))$. Then
$$F'_{\beta,K_x}(x) = \frac{x}{2\beta K_x} - c'_\beta(x) = 0.$$

Since $c'_\beta(w)$ is strictly concave for $w > 0$, we have $c''_\beta(0) > c'_\beta(x)/x$, and, therefore,
$$K_x = \frac{x}{2\beta c'_\beta(x)} > \frac{1}{2\beta c''_\beta(0)} = K_c^{(2)}(\beta).$$

It follows that $x$ is a positive critical point of $F_{\beta,K}$ for $K = K_x > K_c^{(2)}(\beta)$; by the uniqueness of the positive critical point, $x = \tilde{w}(\beta, K_x)$. This completes the proof that the points $\{\tilde{w}(\beta, K), K > K_c^{(2)}(\beta)\}$ span the positive real line. □

PROOF OF THEOREM 3.6. (a) For $0 < K \leq K_c^{(2)}(\beta)$, $F'_{\beta,K}(0) = 0$, and, thus, $w = 0$ is a critical point of $F_{\beta,K}$. We prove that $w = 0$ is the unique global minimum point of $F_{\beta,K}$ by showing that, for $w > 0$, $F'_{\beta,K}(w) > 0$ and for $w < 0$, $F'_{\beta,K}(w) < 0$. Since $c'_{\beta}(w)$ is strictly concave for $w > 0$ [Theorem 3.5(a)], for any $w > 0$, we have $c''_{\beta}(0) > c'_{\beta}(w)/w$. As a result, for all $w > 0$ and all $0 < K \leq K_c^{(2)}(\beta) = 1/(2\beta c''_{\beta}(0))$,

$$F'_{\beta,K}(w) = \frac{w}{2\beta K} - c'_{\beta}(w)$$
$$\geq \frac{w}{2\beta K_c^{(2)}(\beta)} - c'_{\beta}(w) = w c''_{\beta}(0) - c'_{\beta}(w) > 0.$$

On the other hand, since $F'_{\beta,K}$ is an odd function, $F'_{\beta,K}(w) < 0$ for all $w < 0$. Therefore, $w = 0$ is the unique global minimum point of $F_{\beta,K}$. It follows that, for $0 < K \leq K_c^{(2)}(\beta)$, $\tilde{\mathcal{E}}_{\beta,K} = \{0\}$.

(b) For $K > K_c^{(2)}(\beta)$, let $\tilde{w}(\beta, K)$ be the unique positive global minimum point of $F_{\beta,K}$, the existence of which is proved in part (a) of Lemma 3.7, and define $\tilde{z}(\beta, K) = \tilde{w}(\beta, K)/(2\beta K)$. It follows that, for $K > K_c^{(2)}(\beta)$, $\tilde{\mathcal{E}}_{\beta,K} = \{\pm \tilde{z}(\beta, K)\}$.

(c) By part (a) of Lemma 3.7,

$$F'_{\beta,K}(\tilde{w}(\beta, K)) = 0 \quad \text{and} \quad F''_{\beta,K}(\tilde{w}(\beta, K)) > 0.$$

The implicit function theorem implies that, for $K > K_c^{(2)}(\beta)$, $\tilde{w}(\beta, K)$ and, thus, $\tilde{z}(\beta, K)$ are continuously differentiable functions of $K$ and are thus continuous. Straightforward calculations yield

$$\frac{\partial \tilde{w}(\beta, K)}{\partial K} = \frac{\tilde{w}(\beta, K)}{2\beta K^2 F''_{\beta,K}(\tilde{w}(\beta, K))}$$

and

$$\frac{\partial \tilde{z}(\beta, K)}{\partial K} = \frac{2\beta \tilde{w}(\beta, K)}{(2\beta K)^2} \left( \frac{c''_{\beta}(\tilde{w}(\beta, K))}{F''_{\beta,K}(\tilde{w}(\beta, K))} \right).$$

Since $\tilde{w}(\beta, K)$ is positive and both $c''_{\beta}(\tilde{w}(\beta, K)) > 0$ and $F''_{\beta,K}(\tilde{w}(\beta, K)) > 0$, $\tilde{w}(\beta, K)$ and $\tilde{z}(\beta, K)$ are strictly increasing functions for $K > K_c^{(2)}(\beta)$.

As $K \searrow K_c^{(2)}(\beta)$, $\tilde{w}(\beta, K) > 0$, $\tilde{w}(\beta, K)$ is strictly decreasing, and the points $\{\tilde{w}(\beta, K), K > K_c^{(2)}(\beta)\}$ span the positive real line [Lemma 3.7(c)].

We conclude that $\lim_{K \to K_c^{(2)}(\beta)^+} \tilde{w}(\beta, K) = 0$ and thus that $\lim_{K \to K_c^{(2)}(\beta)^+} \tilde{z}(\beta, K) = 0$. This completes the proof of the theorem. $\square$

Theorem 3.6 describes the continuous bifurcation exhibited by $\tilde{\mathcal{E}}_{\beta,K}$ for $0 < \beta \leq \beta_c$. Theorem 3.8 in the next section describes the discontinuous bifurcation exhibited by $\tilde{\mathcal{E}}_{\beta,K}$ for $\beta$ in the complementary region $\beta > \beta_c$.

3.3. *Description of $\tilde{\mathcal{E}}_{\beta,K}$ for $\beta > \beta_c$.* In Theorem 3.2 we state the structure of the set $\mathcal{E}_{\beta,K}$ of canonical equilibrium macrostates for the BEG model with respect to the empirical measures when $\beta > \beta_c$. The main theorem in this subsection, Theorem 3.8, does the same for the set $\tilde{\mathcal{E}}_{\beta,K}$, which has been shown to have the alternative characterization

$$(3.22) \qquad \tilde{\mathcal{E}}_{\beta,K} = \left\{ \frac{w}{2\beta K} \in \mathbb{R} : w \text{ minimizes } F_{\beta,K}(w) \right\}.$$

As in Section 3.2, $F_{\beta,K}(w) = w^2/(4\beta K) - c_\beta(w)$, where $c_\beta$ is defined in (3.4). In Section 3.4 we will prove that there exists a one-to-one correspondence between $\tilde{\mathcal{E}}_{\beta,K}$ and $\mathcal{E}_{\beta,K}$. In Section 3.5 we will use this fact to fully describe the latter set for all $\beta > \beta_c$ and $K > 0$.

Minimum points of $F_{\beta,K}$ satisfy the equation

$$(3.23) \qquad F'_{\beta,K}(w) = \frac{w}{2\beta K} - c'_\beta(w) = 0.$$

In contrast to the previous section, where for $0 < \beta \leq \beta_c$, $c'_\beta(w)$ is strictly concave for $w > 0$, part (b) of Theorem 3.5 states that, for $\beta > \beta_c$, there exists $w_c(\beta) > 0$ such that $c'_\beta(w)$ is strictly convex for $w \in (0, w_c(\beta))$ and strictly concave for $w > w_c(\beta)$. As a result, for $\beta > \beta_c$, we are no longer in the situation of the classical Curie–Weiss model for which the bifurcation with respect to $K$ is continuous. Instead, for $\beta > \beta_c$, as $K$ increases through the critical value $K_c^{(1)}(\beta)$, $\tilde{\mathcal{E}}_{\beta,K}$ exhibits a discontinuous bifurcation.

While the discontinuous bifurcation exhibited by $\tilde{\mathcal{E}}_{\beta,K}$ for $\beta > \beta_c$ is easily observed graphically, the full analytic proof is more complicated than in the case $0 < \beta \leq \beta_c$. As in the previous subsection, we will first motivate this discontinuous bifurcation via a graphical argument. A detailed statement is given in Theorem 3.8.

For $\beta > \beta_c$, we divide the range of the positive parameter $K$ into three intervals separated by the values $K_1 = K_1(\beta)$ and $K_2 = K_2(\beta)$. $K_1$ is defined to be the unique value of $K$ such that the line $w/(2\beta K)$ is tangent to the curve $c'_\beta$ at a point $w_1 = w_1(\beta) > 0$. The existence and uniqueness of $K_1$ and $w_1$ are proved in Lemma 3.9. $K_2$ is defined to be the value of $K$ such that the slopes of the line $w/(2\beta K)$ and the curve $c'_\beta$ at $w = 0$ agree. Specifically,

$$(3.24) \qquad K_2 = \frac{1}{2\beta c''_\beta(0)} = \frac{1}{4\beta e^{-\beta}} + \frac{1}{2\beta}.$$

Figure 2 represents graphically the values of $K_1$ and $K_2$ for $\beta = 4$, showing that $K_1 < K_2$. In Lemma 3.9 it is proved that this inequality holds for all $\beta > \beta_c$.

In each of Figures 3–7, for fixed $\beta > \beta_c$ and for different ranges of values of $K > 0$, the first graph (a) depicts the two components of $F'_{\beta,K}$: the linear component $w/(2\beta K)$ and the nonlinear component $c'_\beta$. The second graph (b) shows the corresponding graph of $F_{\beta,K}$. In these figures the following values of $\beta$ were used: $\beta = 4$ in Figures 3, 5, 6, 7 and $\beta = 2.8$ in Figure 4.

As we see in Figure 3, for $0 < K < K_1$, the linear component intersects the nonlinear component at only the origin and, thus, $F_{\beta,K}$ has a unique global minimum point at $w = 0$. Since $\beta$ is fixed, the graph of the nonlinear component $c'_\beta$ also remains fixed. As $K$ increases, the slope of the linear component $w/(2\beta K)$ decreases, leading to the discontinuous bifurcation in $\mathcal{E}_{\beta,K}$ with respect to $K$.

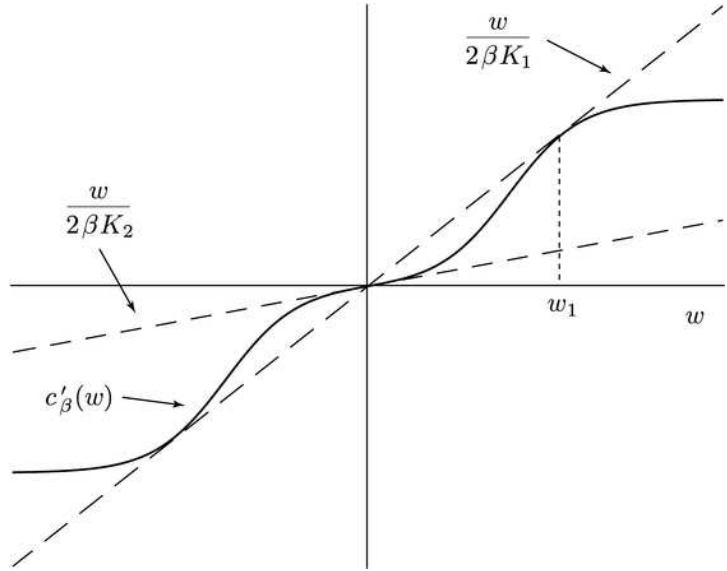

FIG. 2. *Graphical representation of the values $K_1$ and $K_2$ for $\beta = 4$.*

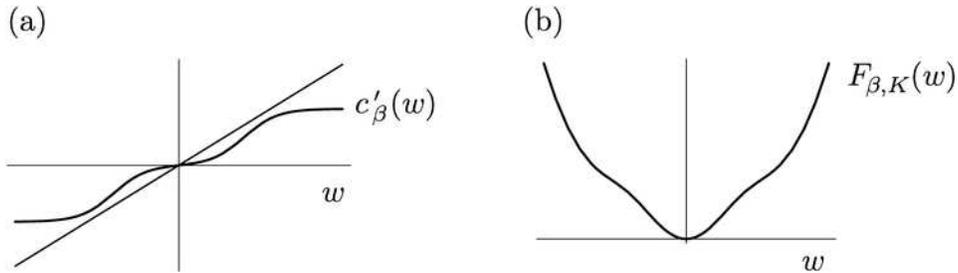

FIG. 3. (a) *Graph of two components of $F'_{\beta,K}$ and* (b) *graph of $F_{\beta,K}$ for $0 < K < K_1$.*

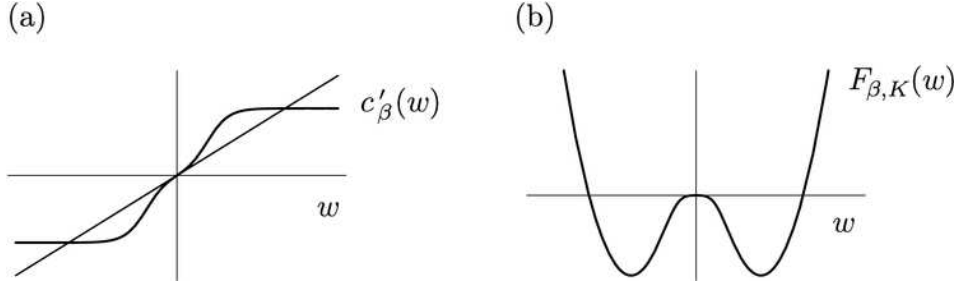

FIG. 4. (a) *Graph of two components of $F'_{\beta,K}$ and* (b) *graph of $F_{\beta,K}$ for $K \geq K_2$.*

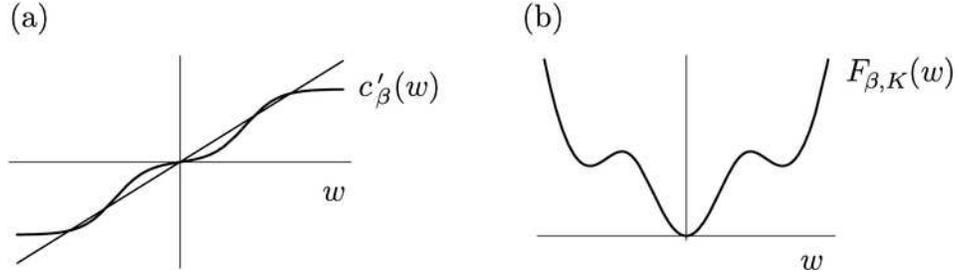

FIG. 5. (a) *Graph of two components of $F'_{\beta,K}$ and* (b) *graph of $F_{\beta,K}$ for $K_1 < K < K_c^{(1)}(\beta)$.*

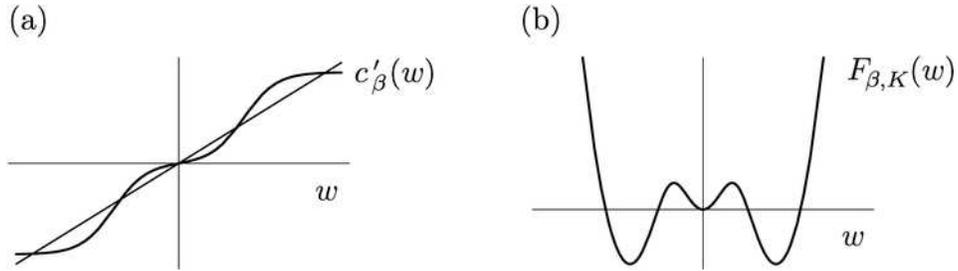

FIG. 6. (a) *Graph of two components of $F'_{\beta,K}$ and* (b) *graph of $F_{\beta,K}$ for $K_c^{(1)}(\beta) < K < K_2$.*

The graph of $F_{\beta,K}$ is depicted in Figure 4 for $K \geq K_2$. We see that $F_{\beta,K}$ has two global minimum points at $w = \pm \tilde{w}(\beta, K)$, where $\tilde{w}(\beta, K)$ is positive. Therefore, for $0 < K \leq K_1$, we have $\tilde{\mathcal{E}}_{\beta,K} = \{0\}$ and for $K \geq K_2$, we have $\tilde{\mathcal{E}}_{\beta,K} = \{\pm \tilde{z}(\beta, K)\}$, where $\tilde{z}(\beta, K) = \tilde{w}(\beta, K)/(2\beta K)$ is positive.

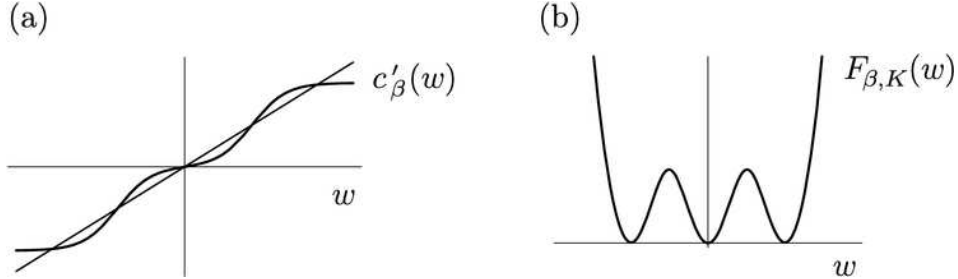

FIG. 7. (a) *Graph of two components of $F'_{\beta,K}$ and* (b) *graph of $F_{\beta,K}$ for $K = K_c^{(1)}(\beta)$.*

Now suppose that $K \in (K_1, K_2)$. In this region there exists $\tilde{w}(\beta, K) > 0$ such that $F_{\beta,K}$ has three local minimum points at $w = 0$ and $w = \pm\tilde{w}(\beta, K)$. As we see in Figure 5, for $K > K_1$ but sufficiently close to $K_1$, $F_{\beta,K}(0) < F_{\beta,K}(\tilde{w}(\beta, K))$; as a result, the unique global minimum point of $F_{\beta,K}$ is $w = 0$. On the other hand, we see in Figure 6 that, for $0 < K < K_2$ but sufficiently close to $K_2$, $F_{\beta,K}(0) > F_{\beta,K}(\tilde{w}(\beta, K))$; as a result, the global minimum points of $F_{\beta,K}$ are $w = \pm\tilde{w}(\beta, K)$. As $K$ increases over the interval $(K_1, K_2)$, $F_{\beta,K}(\tilde{w}(\beta, K))$ decreases continuously (Lemma 3.12). Consequently, as Figure 7 reveals, there exists a critical value $K_c^{(1)}(\beta)$ such that $F_{\beta,K_c^{(1)}(\beta)}(0) = F_{\beta,K_c^{(1)}(\beta)}(\tilde{w}(\beta, K))$; as a result, the global minimum points of $F_{\beta,K_c^{(1)}(\beta)}$ are $w = 0$ and $w = \pm\tilde{w}(\beta, K)$.

We use the same notation $K_c^{(1)}(\beta)$ as for the critical value in Theorem 3.2. As we will later prove, the discontinuous bifurcation in $K$ exhibited by both sets $\mathcal{E}_{\beta,K}$ and $\tilde{\mathcal{E}}_{\beta,K}$ occur at the same point $K_c^{(1)}(\beta)$.

The graphical information just obtained concerning the global minimum points of $F_{\beta,K}$ for $\beta > \beta_c$ motivates the form of $\tilde{\mathcal{E}}_{\beta,K}$ stated in the next theorem. The positive quantity $\tilde{z}(\beta, K)$ equals $\tilde{w}(\beta, K)/(2\beta K)$, where $\tilde{w}(\beta, K)$ is the unique positive global minimum point of $F_{\beta,K}$ for $K \geq K_c^{(1)}(\beta)$ [Lemma 3.10(b)]. According to part (d) of the theorem, $\tilde{z}(\beta, K)$ is a continuous function for $K > K_c^{(1)}(\beta)$, and as $K \to (K_c^{(1)}(\beta))^+$, $\tilde{z}(\beta, K)$ converges to the positive quantity $\tilde{z}(\beta, K_c^{(1)}(\beta))$. Hence, the bifurcation exhibited by $\tilde{\mathcal{E}}_{\beta,K}$ at $K_c^{(1)}(\beta)$ is discontinuous. As we will see in the proof of Theorem 3.8, $K_c^{(1)}(\beta)$ is the unique zero of the function $A(K)$ defined in (3.31) for $K \geq K_1(\beta)$; $K_1(\beta)$ is specified in Lemma 3.9.

THEOREM 3.8. *Define $\tilde{\mathcal{E}}_{\beta,K}$ by* (3.7); *equivalently,*

$$\tilde{\mathcal{E}}_{\beta,K} = \left\{ \frac{w}{2\beta K} \in \mathbb{R} : w \text{ minimizes } F_{\beta,K}(w) \right\}.$$

*For all $\beta > \beta_c = \log 4$, there exists a critical value $K_c^{(1)}(\beta)$ satisfying $K_1 < K_c^{(1)}(\beta) < K_2$ and having the following properties:*

(a) *For $0 < K < K_c^{(1)}(\beta)$, $\tilde{\mathcal{E}}_{\beta,K} = \{0\}$.*
(b) *For $K = K_c^{(1)}(\beta)$, $\tilde{\mathcal{E}}_{\beta,K} = \{0, \pm \tilde{z}(\beta, K)\}$, where $\tilde{z}(\beta, K) > 0$.*
(c) *For $K > K_c^{(1)}(\beta)$, $\tilde{\mathcal{E}}_{\beta,K} = \{\pm \tilde{z}(\beta, K)\}$, where $\tilde{z}(\beta, K) > 0$.*
(d) *For $K \geq K_c^{(1)}(\beta)$, $\tilde{z}(\beta, K)$ is a strictly increasing continuous function, and as $K \to (K_c^{(1)}(\beta))^+$, $\tilde{z}(\beta, K) \to \tilde{z}(\beta, K_c^{(1)}(\beta)) > 0$. Therefore, $\tilde{\mathcal{E}}_{\beta,K}$ exhibits a discontinuous bifurcation at $K_c^{(1)}(\beta)$.*

The proof of the theorem depends on several lemmas. In the first lemma we prove that, for each $\beta > \beta_c$, there exists a unique $K = K_1(\beta)$ such that the line $w/(2\beta K)$ is tangent to the curve $c'_\beta$ at a point $w_1(\beta) > 0$.

LEMMA 3.9.   *For $\beta > \beta_c = \log 4$, we define $c_\beta$ by (3.4), $F_{\beta,K}$ by (3.12), $w_c(\beta)$ by (3.15) and $K_2 = K_2(\beta)$ by (3.24). Then in the set $w > 0$, $K > 0$, there exists a unique solution $(w_1, K_1) = (w_1(\beta), K_1(\beta))$ of*

$$(3.25) \qquad F'_{\beta,K}(w) = \frac{w}{2\beta K} - c'_\beta(w) = 0,$$

$$(3.26) \qquad F''_{\beta,K}(w) = \frac{1}{2\beta K} - c''_\beta(w) = 0.$$

*Furthermore, $w_1 > w_c(\beta)$ and $K_1 < K_2$ for all $\beta > \beta_c$.*

PROOF.   The function $g(w) = wc''_\beta(w) - c'_\beta(w)$ has the properties that $g'(w) = wc'''_\beta(w)$ and that solutions of (3.25)–(3.26) solve $g(w) = 0$. According to part (b) of Theorem 3.5, $c'_\beta(w)$ is strictly convex for $0 < w < w_c(\beta)$ and $c'_\beta(w)$ is strictly concave for $w > w_c(\beta)$; equivalently, $c'''_\beta(w) > 0$ for $0 < w < w_c(\beta)$ and $c'''_\beta(w) < 0$ for $w > w_c(\beta)$. Therefore,

$$(3.27) \quad g'(w) > 0 \quad \text{for } 0 < w < w_c(\beta) \quad \text{and} \quad g'(w) < 0 \quad \text{for } w > w_c(\beta).$$

Since

$$wc''_\beta(w) = \frac{2we^{-\beta} \cosh w + 4we^{-2\beta}}{[1 + 2e^{-\beta} \cosh w]^2},$$

we see that $\lim_{w \to \infty} wc''_\beta(w) = 0$. It follows that

$$\lim_{w \to \infty} g(w) = \lim_{w \to \infty} wc''_\beta(w) - \lim_{w \to \infty} c'_\beta(w) = -1.$$

This limit and the fact that $g(0) = 0$, combined with the continuity of $g$ and (3.27), imply that there exists a unique $w_1 > w_c(\beta)$ such that $g(w_1) = 0$;

that is,

$$\frac{c'_\beta(w_1)}{w_1} = c''_\beta(w_1). \tag{3.28}$$

Substituting $w_1$ into (3.25) and (3.26), we define

$$K_1 = \frac{1}{2\beta c''_\beta(w_1)} = \frac{w_1}{2\beta c'_\beta(w_1)}. \tag{3.29}$$

The pair $(w_1, K_1)$ is a solution of (3.25)–(3.26) in the set $w > 0, K > 0$. If $(\hat{w}, \hat{K})$ is another solution of (3.25)–(3.26) in this set, then $\hat{w}$ solves $g(w) = 0$, a contradiction to the fact that $w_1$ is the unique positive solution of $g(w) = 0$. It follows that $(w_1, K_1)$ is the unique solution of (3.25)–(3.26) in the set $w > 0, K > 0$.

We complete the proof by showing that $K_1 < K_2$. Since $K_2 = 1/(2\beta c''_\beta(0))$, we are done if we show that $c''_\beta(w_1) > c''_\beta(0)$. By the mean value theorem and (3.28), there exists $\alpha \in (0, w_1)$ such that

$$c''_\beta(\alpha) = \frac{c'_\beta(w_1)}{w_1} = c''_\beta(w_1). \tag{3.30}$$

We claim that $\alpha < w_c(\beta)$. If $\alpha \geq w_c(\beta)$, then since $c'_\beta$ is strictly concave on $(w_c(\beta), \infty)$, the inequalities $w_c(\beta) \leq \alpha < w_1$ imply that $c''_\beta(\alpha) > c''_\beta(w_1)$. Because this contradicts (3.30), we conclude that $\alpha < w_c(\beta)$. This inequality in combination with the strict convexity of $c'_\beta$ on $(0, w_c(\beta))$ and (3.30) yields $c''_\beta(0) < c''_\beta(\alpha) = c''_\beta(w_1)$. The proof of the lemma is complete. $\square$

We next state two lemmas that are analogous to Lemma 3.7 and part (c) of Theorem 3.6. Before stating them, we need some preliminaries. In Lemma 3.9, we proved that, for $\beta > \beta_c$, equations (3.25)–(3.26) have a unique solution $(w_1, K_1) = (w_1(\beta), K_1(\beta))$ in the set $w > 0, K > 0$ and that $w_1 > w_c(\beta)$; according to (3.25),

$$F''_{\beta, K_1}(w_1) = \frac{1}{2\beta K_1} - c''_\beta(w_1) = 0.$$

In addition, for $0 < \beta \leq \beta_c$, the quantity $K_c^{(2)}(\beta) = 1/(2\beta c''_\beta(0))$ introduced in (3.19) has the property that

$$F''_{\beta, K_c^{(2)}(\beta)}(0) = \frac{1}{2\beta K_c^{(2)}(\beta)} - c''_\beta(0) = 0.$$

For $\beta > \beta_c$, $c'_\beta(w)$ is strictly concave for $w > w_c(\beta)$ [Theorem 3.5(b)]. Thus for $w \geq w_1$, the graph of $F_{\beta, K_1}(w)$ over the interval $[w_1, \infty)$ for $\beta > \beta_c$ (Figure 2) is similar to that of $F_{\beta, K_c^{(2)}(\beta)}(w)$ over the interval $[0, \infty)$ for $0 < \beta \leq \beta_c$

[Figure 1(b)]. Specifically, for $\beta > \beta_c$ and $w \in [w_1, \infty)$, the graph of

$$F_{\beta, K_1}(w) = \int_{w_1}^{w} \left( \frac{x}{2\beta K_1} - c'_\beta(x) \right) dx + F_{\beta, K_1}(w_1)$$

is determined by the difference between the strictly concave function $c'_\beta(w)$ and the linear function $w/(2\beta K_1)$, which is tangent to $c'_\beta$ at $w = w_1$. Similarly, for $0 < \beta \leq \beta_c$ and $w \in [0, \infty)$, the graph of

$$F_{\beta, K_c^{(2)}(\beta)}(w) = \int_0^w \left( \frac{x}{2\beta K_c^{(2)}(\beta)} - c'_\beta(x) \right) dx$$

is determined by the difference between the strictly concave function $c'_\beta(w)$ [Theorem 3.5(a)] and the linear function $w/(2\beta K_c^{(2)}(\beta))$, which is tangent to $c'_\beta$ at $w = 0$.

As we saw in Section 3.2 for $0 < \beta \leq \beta_c$, as $K$ increases from $K_c^{(2)}(\beta)$ and thus the slope of the line $w/(2\beta K)$ decreases, $F_{\beta, K}$ develops a unique positive local minimum point $\tilde{w}(\beta, K)$, which is shown to be the unique global minimum point on the interval $[0, \infty)$ [Lemma 3.7(b)]. This can be seen graphically in Figure 1(c). Similarly, as Figures 4(b)–7(b) illustrate, for $\beta > \beta_c$, as $K$ increases from $K_1$, $F_{\beta, K}$ develops a unique positive local minimum point $\tilde{w}(\beta, K)$. As in part (b) of Lemma 3.7, $\tilde{w}(\beta, K)$ can be shown to be the unique global minimum point on the interval $[w_1, \infty)$. However, it is not a global minimum point on the entire halfline $[0, \infty)$ unless $F_{\beta, K}(\tilde{w}(\beta, K)) \leq 0 = F_{\beta, K}(0)$; in fact, this inequality is valid only for all $K$ sufficiently large. When $F_{\beta, K}(\tilde{w}(\beta, K)) > 0 = F_{\beta, K}(0)$, which holds for all $K > K_1$ sufficiently close to $K_1$, 0 is the unique global minimum point of $F_{\beta, K}$.

Because the behavior of the function $F_{\beta, K}$ over the interval $[w_1, \infty)$ for $\beta > \beta_c$ is similar to that of $F_{\beta, K}$ over the interval $[0, \infty)$ for $0 < \beta \leq \beta_c$, the proofs of Lemma 3.10 and Lemma 3.11 are analogous, respectively, to the proofs of Lemma 3.7 and part (c) of Theorem 3.6. Therefore, we state these new lemmas without proof.

In Lemma 3.10 we state the existence and two properties of a positive critical point $\tilde{w}(\beta, K)$ of $F_{\beta, K}$ for each $K > K_1$.

LEMMA 3.10. *For $\beta > \beta_c = \log 4$, define $F_{\beta, K}$ by* (3.12) *and let* $(w_1, K_1) = (w_1(\beta), K_1(\beta))$ *be the unique solution of* (3.25)–(3.26) *in the set $w > 0, K > 0$ (Lemma* 3.9*). The following conclusions hold:*

(a) *For each $K > K_1$, $F_{\beta, K}$ has a critical point $\tilde{w}(\beta, K) > w_1$ satisfying*

$$F'_{\beta, K}(\tilde{w}(\beta, K)) = 0 \quad and \quad F''_{\beta, K}(\tilde{w}(\beta, K)) > 0.$$

(b) *For each $K > K_1$, $F_{\beta, K}$ has unique nonzero local minimum points at $w = \pm \tilde{w}(\beta, K)$.*

(c) *The points $\{\tilde{w}(\beta, K), K > K_1\}$ span the interval $(w_1, \infty)$; that is, for each $x > w_1$, there exists $K > K_1$ such that $x = \tilde{w}(\beta, K)$.*

The next lemma states continuity and related properties of $\tilde{w}(\beta, K)$ and $\tilde{z}(\beta, K)$ that are similar to properties of the analogous quantities for $0 < \beta \leq \beta_c$ [Theorem 3.6(c)].

LEMMA 3.11. *For $\beta > \beta_c = \log 4$ and $K > K_1$, let $\tilde{w}(\beta, K)$ be the unique positive local minimum point of $F_{\beta, K}$ considered in Lemma 3.10. Then for $K > K_1$, $\tilde{w}(\beta, K)$ and $\tilde{z}(\beta, K) = \tilde{w}(\beta, K)/(2\beta K)$ are continuous, strictly increasing functions of $K$ and $\lim_{K \to K_1^+} \tilde{w}(\beta, K) = w_1$.*

We fix $\beta > \beta_c$. The proof of Theorem 3.8 also makes use of the function

$$(3.31) \qquad D(K) = \begin{cases} F_{\beta, K_1}(w_1), & \text{if } K = K_1, \\ F_{\beta, K}(\tilde{w}(\beta, K)), & \text{if } K > K_1. \end{cases}$$

The quantity $\tilde{w}(\beta, K)$ is the unique positive local minimum point of $F_{\beta, K}$, the existence of which is given in Lemma 3.10.

LEMMA 3.12. *For $\beta > \beta_c = \log 4$, the function $D(K)$ defined in (3.31) is continuous and strictly decreasing on its domain $[K_1(\beta), \infty)$.*

PROOF. Since $F_{\beta, K}(w)$ is a continuous function of $w$ and $\tilde{w}(\beta, K)$ is a continuous function of $K$ (Lemma 3.11), $D(K)$ is continuous for $K > K_1$. Furthermore, by part (c) of Lemma 3.11, $\lim_{K \to K_1^+} \tilde{w}(\beta, K) = w_1$ and, thus, $\lim_{K \to K_1^+} F_{\beta, K}(\tilde{w}(\beta, K)) = F_{\beta, K_1}(w_1)$. We conclude that $D(K)$ is continuous on $[K_1, \infty)$.

We now prove that $D(K)$ is strictly decreasing on $[K_1, \infty)$. For $K > K_1$, we have

$$\frac{\partial F_{\beta, K}}{\partial w}(\tilde{w}(\beta, K)) = 0$$

by part (a) of Lemma 3.10. As in the proof of part (c) of Theorem 3.6, one can show that $\tilde{w}(\beta, K)$ is continuously differentiable for $K > K_1$. Hence, for $K > K_1$,

$$\begin{aligned} D'(K) &= \frac{dF_{\beta, K}(\tilde{w}(\beta, K))}{dK} \\ &= \frac{\partial F_{\beta, K}}{\partial K}(\tilde{w}(\beta, K)) + \frac{\partial F_{\beta, K}}{\partial w}(\tilde{w}(\beta, K)) \cdot \frac{\partial \tilde{w}(\beta, K)}{\partial K} \\ &= -\frac{[\tilde{w}(\beta, K)]^2}{4\beta K^2} < 0. \end{aligned}$$

This completes the proof. □

PROOF OF THEOREM 3.8. As we showed in Lemma 3.9, for $\beta > \beta_c$, equations (3.25)–(3.26) have a unique solution $(w_1, K_1) = (w_1(\beta), K_1(\beta))$ in the set $w > 0, K > 0$. In addition, $K_1 < K_2 = 1/(2\beta c''_\beta(0))$. We start the proof of Theorem 3.8 by proving the following two facts:

1. For $0 < K \leq K_1$, $F_{\beta,K}$ has a unique global minimum point at $w = 0$ [Figures 2 and 3(b)].
2. For $K \geq K_2$, $F_{\beta,K}$ has unique global minimum points at $w = \pm\tilde{w}(\beta, K)$ (Figure 7).

According to Lemma 3.10, $F'_{\beta,K}(w_1) = 0$. Using concavity properties of $c'_\beta(w)$ established in part (b) of Theorem 3.5 and calculations similar to those used to establish other results in this and the preceding section, one shows that, for $0 < K < K_1$, $F'_{\beta,K}(w) > 0$ for all $w > 0$ and that $F'_{\beta,K_1}(w) > 0$ for all $w > 0$, $w \neq w_1$. These properties, which can be seen in Figure 2 and Figure 3(a), are proved in detail in Lemma 2.3.10 in [21]. By symmetry, for $0 < K < K_1$, $F'_{\beta,K}(w) < 0$ for all $w < 0$ and $F'_{\beta,K}(w) < 0$ for all $w < 0, w \neq -w_1$. It follows that, for $0 < K \leq K_1$, $F_{\beta,K}$ is strictly decreasing for $w < 0$ and strictly increasing for $w > 0$. We conclude that, for $0 < K \leq K_1$, $F_{\beta,K}$ has a unique global minimum point at $w = 0$, as claimed in fact 1.

Since $\lim_{|w|\to\infty} F_{\beta,K}(w) = \infty$, the global minimum values of $F_{\beta,K}$ must be attained at local minimum points of the function. Lemma 3.10 states that, for $K > K_1$, $w = \pm\tilde{w}(\beta, K)$ are the unique nonzero local minimum points of $F_{\beta,K}$. Therefore, we prove that, for $K \geq K_2$, $F_{\beta,K}$ has unique global minimum points at $w = \pm\tilde{w}(\beta, K)$ by proving that $w = 0$ is a local maximum point of $F_{\beta,K}$. According to part (b) of Theorem 3.5, $c'_\beta(w)$ is strictly convex for $0 < w < w_c(\beta)$. Therefore, for $K \geq K_2$ and $w \in (0, w_c(\beta))$,

$$F'_{\beta,K}(w) = \frac{w}{2\beta K} - c'_\beta(w) = \frac{w}{2\beta K} - \int_0^w c''_\beta(x)\,dx$$
$$< \frac{w}{2\beta K} - c''_\beta(0)w = \frac{w}{2\beta K} - \frac{w}{2\beta K_2} \leq 0;$$

that is, for $w \in (0, w_c(\beta))$, $F'_{\beta,K}(w) < 0$. By symmetry, for $w \in (-w_c(\beta), 0)$, $F'_{\beta,K}(w) > 0$. It follows that $w = 0$ is a local maximum point of $F_{\beta,K}$. Therefore, as claimed in fact 2, for $K \geq K_2$, $F_{\beta,K}$ has unique global minimum points at $w = \pm\tilde{w}(\beta, K)$.

For $K_1 < K < K_2$, $F_{\beta,K}$ has exactly three local minimum points at $w = 0$ and $w = \pm\tilde{w}(\beta, K)$. Since global minimum values of $F_{\beta,K}$ must be attained at local minimum points of the function and $F_{\beta,K}$ is symmetric, finding global minimum points of $F_{\beta,K}$ requires comparing the values of the function at $w = 0$ and at $w = \tilde{w}(\beta, K)$.

Since for $0 < K \leq K_1$ $F_{\beta,K}$ has a unique global minimum point at $w = 0$, we have

$$D(K_1) = F_{\beta,K_1}(w_1) > F_{\beta,K_1}(0) = 0.$$

Similarly, since for $K \geq K_2$ $F_{\beta,K}$ has unique global minimum points at $w = \pm\tilde{w}(\beta, K)$, for $K \geq K_2$, we have

$$D(K) = F_{\beta,K}(\tilde{w}(\beta, K)) < F_{\beta,K}(0) = 0.$$

Since $D(K)$ is continuous and strictly decreasing for $K > K_1$ (Lemma 3.12), there exists a unique critical value $K_c^{(1)}(\beta)$ satisfying $K_1 < K_c^{(1)}(\beta) < K_2$ and having the following properties:

(i) For $K_1 < K < K_c^{(1)}(\beta)$,

$$F_{\beta,K}(\tilde{w}(\beta, K)) = D(K) > 0 = F_{\beta,K_1}(0),$$

and, thus, $F_{\beta,K}$ has a unique global minimum point at $w = 0$.

(ii) For $K = K_c^{(1)}(\beta)$,

$$F_{\beta,K}(\tilde{w}(\beta, K)) = D(K) = 0 = F_{\beta,K_1}(0),$$

and, thus, $F_{\beta,K}$ has three global minimum points at $w = 0, \pm\tilde{w}(\beta, K)$.

(iii) For $K_c^{(1)}(\beta) < K < K_2$,

$$F_{\beta,K}(\tilde{w}(\beta, K)) = D(K) < 0 = F_{\beta,K_1}(0),$$

and, thus, $F_{\beta,K}$ has two global minimum points at $w = \pm\tilde{w}(\beta, K)$.

We define $\tilde{z}(\beta, K) = \tilde{w}(\beta, K)/(2\beta K)$. The form of $\tilde{\mathcal{E}}_{\beta,K}$ given in parts (a)–(c) of Theorem 3.8 follows from the information on the global minimum points of $F_{\beta,K}$ just given in items (i)–(iii) and from facts 1 and 2 stated at the start of the proof. In addition, the positivity of $\tilde{z}(\beta, K_c^{(1)}(\beta))$ is a consequence of the positivity of $\tilde{w}(\beta, K_c^{(1)}(\beta))$. Since by Lemma 3.11 $\tilde{z}(\beta, K)$ is a strictly increasing function for $K \geq K_c^{(1)}(\beta)$, part (d) of the theorem is also proved. The proof of Theorem 3.8 is now complete. $\square$

Together, Theorems 3.6 and 3.8 give a full description of the set $\tilde{\mathcal{E}}_{\beta,K}$ for all values of $\beta$ and $K$. In the next section, we use the contraction principle to lift our results concerning the structure of the set $\tilde{\mathcal{E}}_{\beta,K}$ up to the level of the empirical measures, making use of a one-to-one correspondence between the points in the two sets $\tilde{\mathcal{E}}_{\beta,K}$ and $\mathcal{E}_{\beta,K}$ of canonical equilibrium macrostates.

3.4. *One-to-one correspondence between $\mathcal{E}_{\beta,K}$ and $\tilde{\mathcal{E}}_{\beta,K}$.* We start by recalling the definitions of the sets $\mathcal{E}_{\beta,K}$ and $\tilde{\mathcal{E}}_{\beta,K}$:

$$\mathcal{E}_{\beta,K} = \{\nu \in \mathcal{P}(\Lambda) : \nu \text{ minimizes } R(\nu|\rho) + \beta f_K(\nu)\} \tag{3.32}$$

and

$$\tilde{\mathcal{E}}_{\beta,K} = \{z \in [-1,1] : z \text{ minimizes } J_\beta(z) - \beta K z^2\}. \tag{3.33}$$

In the definition of $\mathcal{E}_{\beta,K}$, $R(\mu|\rho)$ is the relative entropy of $\mu$ with respect to $\rho = \frac{1}{3}(\delta_{-1} + \delta_0 + \delta_1)$ and $f_K(\mu)$ is the function defined in (2.3). In the definition of $\tilde{\mathcal{E}}_{\beta,K}$, $J_\beta$ is the Cramér rate function defined in (3.3). We now state the one-to-one correspondence between the points in $\tilde{\mathcal{E}}_{\beta,K}$ and the points in $\mathcal{E}_{\beta,K}$. According to Theorems 3.6 and 3.8, $\tilde{\mathcal{E}}_{\beta,K}$ consists of either $1, 2$ or $3$ points.

THEOREM 3.13. *Fix $\beta > 0$ and $K > 0$ and suppose that $\tilde{\mathcal{E}}_{\beta,K} = \{z_\alpha\}_{\alpha=1}^r$, $r = 1, 2$ or $3$. Define $\nu_\alpha, \alpha = 1, \ldots, r$, to be measures in $\mathcal{P}(\Lambda)$ with densities*

$$\frac{d\nu_\alpha}{d\rho_\beta}(y) = \exp(t_\alpha y) \cdot \frac{1}{\int_\Lambda \exp(t_\alpha y) \rho_\beta(dy)}, \tag{3.34}$$

*where $t_\alpha$ is chosen such that $\int_\Lambda y \nu_\alpha(dy) = z_\alpha$. Then for each $\alpha = 1, \ldots, r$, $t_\alpha$ exists and is unique, and $\mathcal{E}_{\beta,K}$ consists of the unique elements $\nu_\alpha, \alpha = 1, \ldots, r$. Furthermore, $t_\alpha = 2\beta K z_\alpha$ for $\alpha = 1, \ldots, r$.*

For $z \in [-1, 1]$, we define

$$A(z) = \left\{\mu \in \mathcal{P}(\Lambda) : \int_\Lambda y\mu(dy) = z\right\}. \tag{3.35}$$

The proof of the theorem depends on the following two lemmas. Both lemmas use the contraction principle ([9], Theorem VIII.3.1), which states that, for all $z \in [-1, 1]$,

$$J_\beta(z) = \min\{R(\mu|\rho_\beta) : \mu \in A(z)\}. \tag{3.36}$$

LEMMA 3.14. *For $\beta > 0$ and $K > 0$,*

$$\min_{\mu \in \mathcal{P}(\Lambda)} \left\{R(\mu|\rho_\beta) - \beta K \left(\int_\Lambda y\mu(dy)\right)^2\right\} = \min_{|z| \leq 1}\{J_\beta(z) - \beta K z^2\}.$$

PROOF. The contraction principle (3.36) implies that

$$\min_{\mu \in \mathcal{P}(\Lambda)} \left\{ R(\mu|\rho_\beta) - \beta K \left( \int_\Lambda y\mu(dy) \right)^2 \right\}$$

$$= \min_{|z| \leq 1} \left( \min \left\{ R(\mu|\rho_\beta) - \beta K \left( \int_\Lambda y\mu(dy) \right)^2 : \mu \in A(z) \right\} \right)$$

$$= \min_{|z| \leq 1} (\min\{R(\mu|\rho_\beta) : \mu \in A(z)\} - \beta K z^2)$$

$$= \min_{|z| \leq 1} \{J_\beta(z) - \beta K z^2\}.$$

This completes the proof. □

The second lemma shows that the mean of any measure $\nu \in \mathcal{E}_{\beta,K}$ is an element of $\tilde{\mathcal{E}}_{\beta,K}$.

LEMMA 3.15. *Fix $\beta > 0$ and $K > 0$. Given $\nu \in \mathcal{E}_{\beta,K}$, we define $\tilde{z} = \int_\Lambda y\nu(dy)$, where $\Lambda = \{-1, 0, 1\}$. Then $\tilde{z} \in \tilde{\mathcal{E}}_{\beta,K}$.*

PROOF. Since $\nu \in \mathcal{E}_{\beta,K}$, $\nu$ is a global minimum point of $R(\mu|\rho_\beta) - \beta K(\int_\Lambda y\mu(dy))^2$. Thus, for all $\mu \in \mathcal{P}(\Lambda)$,

$$R(\nu|\rho_\beta) - \beta K \left( \int_\Lambda y\nu(dy) \right)^2$$

$$= R(\nu|\rho_\beta) - \beta K \tilde{z}^2 \leq R(\mu|\rho_\beta) - \beta K \left( \int_\Lambda y\mu(dy) \right)^2.$$

In particular, this inequality holds for any $\mu$ that satisfies $\int_\Lambda y\mu(dy) = \tilde{z}$. For such $\mu$, the last display becomes

$$R(\nu|\rho_\beta) \leq R(\mu|\rho_\beta).$$

Thus, $\nu$ satisfies

$$R(\nu|\rho_\beta) = \min\{R(\mu|\rho_\beta) : \mu \in A(\tilde{z})\},$$

where $A(\tilde{z})$ is defined in (3.35). The contraction principle (3.36) and Lemma 3.14 imply that

$$J_\beta(\tilde{z}) - \beta K \tilde{z}^2 = R(\nu|\rho_\beta) - \beta K \left( \int_\Lambda y\nu(dy) \right)^2$$

$$= \min_{\mu \in \mathcal{P}(\Lambda)} \left\{ R(\mu|\rho_\beta) - \beta K \left( \int_\Lambda y\mu(dy) \right)^2 \right\}$$

$$= \min_{|z| \leq 1} \{J_\beta(z) - \beta K z^2\}.$$

Therefore, $\tilde{z} \in \tilde{\mathcal{E}}_{\beta,K}$, as claimed. This completes the proof. $\square$

We next prove Theorem 3.13.

PROOF OF THEOREM 3.13. A short calculation shows that, for any $\mu \in \mathcal{P}(\Lambda)$,

$$R(\mu|\rho) + \beta f_K(\mu) - \inf_{\nu \in \mathcal{P}(\Lambda)} \{R(\nu|\rho) + \beta f_K(\nu)\}$$

$$= R(\mu|\rho_\beta) - \beta K \left(\int_\Lambda y\mu(dy)\right)^2 - \inf_{\nu \in \mathcal{P}(\Lambda)} \left\{R(\nu|\rho_\beta) - \beta K \left(\int_\Lambda y\nu(dy)\right)^2\right\}.$$

Hence, we obtain the following alternate characterization of $\mathcal{E}_{\beta,K}$:

$$(3.37) \quad \mathcal{E}_{\beta,K} = \left\{\nu \in \mathcal{P}(\Lambda) : \nu \text{ minimizes } R(\nu|\rho_\beta) - \beta K \left(\int_\Lambda y\nu(dy)\right)^2\right\}.$$

We first show for each $\alpha = 1, \ldots, r$ and $z_\alpha \in \tilde{\mathcal{E}}_{\beta,K}$, $\nu_\alpha$ is the unique global minimum point of $R(\mu|\rho_\beta) - \beta K (\int_\Lambda y\mu(dy))^2$ over

$$A(z_\alpha) = \left\{\mu \in \mathcal{P}(\Lambda) : \int_\Lambda y\mu(dy) = z_\alpha\right\}.$$

We then prove that

$$\inf_{\mu \in A(z_\alpha)} \left\{R(\mu|\rho_\beta) - \beta K \left(\int_\Lambda y\mu(dy)\right)^2\right\}$$

$$= \inf_{\mu \in A(z_\ell)} \left\{R(\mu|\rho_\beta) - \beta K \left(\int_\Lambda y\mu(dy)\right)^2\right\}$$

for all $\alpha, \ell = 1, \ldots, r$. It will then follow that $\{\nu_\alpha\}_{\alpha=1}^r$ equals the set of global minimum points of $R(\mu|\rho_\beta) - \beta K(\int_\Lambda y\mu(dy))^2$ over the set $A = \bigcup_{\alpha=1}^r A(z_\alpha)$. Finally, by showing that all the global minimum points of $R(\mu|\rho_\beta) - \beta K(\int_\Lambda y\mu(dy))^2$ lie in $A$, we will complete the proof that $\mathcal{E}_{\beta,K} = \{\nu_\alpha\}_{\alpha=1}^r$. If $r = 2$ or 3, then since $\int_\Lambda y\nu_\alpha(dy) = z_\alpha$, it is clear that if $z_\alpha \neq z_\ell$, then $\nu_\alpha \neq \nu_\ell$.

By Theorem VIII.3.1 in [9], for each $\alpha = 1, \ldots, r$, the point $t_\alpha$ in the statement of Theorem 3.13 exists and is unique,

$$(3.38) \qquad\qquad\qquad J_\beta(z_\alpha) = R(\nu_\alpha|\rho_\beta),$$

and $R(\mu|\rho_\beta)$ attains its infimum over $A(z_\alpha)$ at the unique measure $\nu_\alpha$. Therefore, for each $\alpha = 1, \ldots, r$, $\nu_\alpha$ is the unique global minimum point of $R(\mu|\rho_\beta) - \beta K(\int_\Lambda y\mu(dy))^2$ over $A(z_\alpha)$.

We next show that

$$\inf_{\mu \in A(z_\alpha)} \left\{R(\mu|\rho_\beta) - \beta K \left(\int_\Lambda y\mu(dy)\right)^2\right\}$$

$$= \inf_{\mu \in A(z_\ell)} \left\{R(\mu|\rho_\beta) - \beta K \left(\int_\Lambda y\mu(dy)\right)^2\right\}$$

for all $\alpha, \ell = 1, \ldots, r$. Since $z_\alpha, z_\ell \in \tilde{\mathcal{E}}_{\beta,K}$, $z_\alpha$ and $z_\ell$ are global minimum points of $J_\beta(z) - \beta K z^2$. Thus, by (3.38), we have

$$\inf_{\mu \in A(z_\alpha)} \left\{ R(\mu|\rho_\beta) - \beta K \left( \int_\Lambda y \mu(dy) \right)^2 \right\}$$
$$= \inf_{\mu \in A(z_\alpha)} R(\mu|\rho_\beta) - \beta K z_\alpha^2$$
$$= J_\beta(z_\alpha) - \beta K z_\alpha^2$$
$$= J_\beta(z_\ell) - \beta K z_\ell^2$$
$$= \inf_{\mu \in A(z_\ell)} R(\mu|\rho_\beta) - \beta K z_\ell^2$$
$$= \inf_{\mu \in A(z_\ell)} \left\{ R(\mu|\rho_\beta) - \beta K \left( \int_\Lambda y \mu(dy) \right)^2 \right\}.$$

As a result, $\{\nu_\alpha\}_{\alpha=1}^r$ equals the set of global minimum points of $R(\mu|\rho_\beta) - \beta K (\int_\Lambda y \mu(dy))^2$ over the set $A = \bigcup_{\alpha=1}^r A(z_\alpha)$.

Last, we show $R(\mu|\rho_\beta) - \beta K (\int_\Lambda y \mu(dy))^2$ attains its global minimum at points in $A$. Let $\sigma$ be a global minimum point of $R(\mu|\rho_\beta) - \beta K (\int_\Lambda y \mu(dy))^2$. By (3.37), this implies that $\sigma \in \mathcal{E}_{\beta,K}$. Define $\zeta = \int_\Lambda y \sigma(dy)$. Then Lemma 3.15 implies that $\zeta \in \tilde{\mathcal{E}}_{\beta,K}$ and, thus, that $\zeta = z_\alpha$ for some $\alpha = 1, \ldots, r$. It follows that $\sigma \in A(z_\alpha) \subset A$ for some $\alpha = 1, \ldots, r$.

The last step is to prove that $t_\alpha = 2\beta K z_\alpha$ for $\alpha = 1, \ldots, r$. From definition (3.4), we have

$$c'_\beta(t_\alpha) = \int_\Lambda y \nu_\alpha(dy) = z_\alpha.$$

In turn, the inverse relationship (3.5) implies that

$$t_\alpha = (c'_\beta)^{-1}(z_\alpha) = J'_\beta(z_\alpha).$$

Therefore, since $z_\alpha \in \tilde{\mathcal{E}}_{\beta,K}$, the definition (3.33) guarantees that $z_\alpha$ is a critical point of $J_\beta(z) - \beta K z^2$. Thus,

(3.39) $$t_\alpha = J'_\beta(z_\alpha) = 2\beta K z_\alpha.$$

This completes the proof of Theorem 3.13. □

In the next section we use Theorem 3.13 to prove Theorems 3.1 and 3.2.

3.5. *Proofs of Theorems 3.1 and 3.2.* Theorem 3.1 gives the structure of the set $\mathcal{E}_{\beta,K}$ of canonical equilibrium macrostates, pointing out the continuous bifurcation exhibited by that set for $0 < \beta \leq \beta_c = \log 4$. The structure of $\mathcal{E}_{\beta,K}$ for $\beta > \beta_c$, given in Theorem 3.2, features a discontinuous bifurcation in $K$. The proofs of these theorems are immediate from Theorems 3.6 and 3.8, respectively, which give the structure of $\tilde{\mathcal{E}}_{\beta,K}$ for $0 < \beta \leq \beta_c$ and for $\beta > \beta_c$, and from Theorem 3.13, which states a one-to-one correspondence between $\tilde{\mathcal{E}}_{\beta,K}$ and $\mathcal{E}_{\beta,K}$.

Before proving Theorems 3.1 and 3.2, it is useful to express the measures $\rho_\beta$ and $\nu_\alpha$ in Theorem 3.13 in the forms $\rho_\beta = \rho_{\beta,-1}\delta_{-1} + \rho_{\beta,0}\delta_0 + \rho_{\beta,1}\delta_1$ and $\nu_\alpha = \nu_{\alpha,-1}\delta_{-1} + \nu_{\alpha,0}\delta_0 + \nu_{\alpha,1}\delta_1$, respectively. Since $t_\alpha = 2\beta K z_\alpha$, in terms of $z_\alpha \in \tilde{\mathcal{E}}_{\beta,K}$ we have

$$\rho_{\beta,-1} = \frac{e^{-\beta}}{1+2e^{-\beta}}, \qquad \rho_{\beta,0} = \frac{1}{1+2e^{-\beta}}, \qquad \rho_{\beta,1} = \frac{e^{-\beta}}{1+2e^{-\beta}}$$

and

$$\nu_{\alpha,-1} = \frac{e^{-2\beta K z_\alpha - \beta}}{C(\beta,K)}, \qquad \nu_{\alpha,0} = \frac{1}{C(\beta,K)}, \qquad \nu_{\alpha,1} = \frac{e^{2\beta K z_\alpha - \beta}}{C(\beta,K)}.$$

Here

$$C(\beta,K) = e^{-2\beta K z_\alpha - \beta} + e^{2\beta K z_\alpha - \beta} + 1.$$

In particular, $\nu_\alpha = \rho_\beta$ when $z_\alpha = 0$.

We first indicate how Theorem 3.1 follows from Theorem 3.6. Fix $0 < \beta \leq \beta_c$. The critical value $K_c^{(2)}(\beta)$ in Theorem 3.1 coincides with the value $K_c^{(2)}(\beta)$ in Theorem 3.6. For $0 < K \leq K_c^{(2)}(\beta)$, part (a) of Theorem 3.6 indicates that $\tilde{\mathcal{E}}_{\beta,K} = \{0\}$; hence, $\mathcal{E}_{\beta,K} = \{\rho_\beta\}$. For $K > K_c^{(2)}(\beta)$, part (b) of Theorem 3.6 indicates that $\tilde{\mathcal{E}}_{\beta,K} = \{\pm\tilde{z}(\beta,K)\}$, where $\tilde{z}(\beta,K) > 0$. It follows that the measures $\nu^+(\beta,K)$ and $\nu^-(\beta,K)$ in part (a)(ii) of Theorem 3.1 are given by (3.34) with $z_\alpha = \tilde{z}(\beta,K)$ and $z_\alpha = -\tilde{z}(\beta,K)$, respectively. Since $\tilde{z}(\beta,K) > 0$, it follows that $\nu^+(\beta,K) \neq \nu^-(\beta,K) \neq \rho_\beta$. Finally, part (c) of Theorem 3.6 allows us to conclude that, for each choice of sign, $\nu^\pm(\beta,K)$ is a continuous functions for $K > K_c^{(2)}(\beta)$ and that as $K \to (K_c^{(2)}(\beta))^+$, $\nu^\pm \to \rho_\beta$. This completes the proof of Theorem 3.1.

In a completely analogous way, Theorem 3.2, including the discontinuous bifurcation noted in part (c) of the theorem, follows from Theorem 3.8.

In this section we have completely analyzed the structure of the set $\mathcal{E}_{\beta,K}$ of canonical equilibrium macrostates. In particular, we discovered that, for $0 < \beta \leq \beta_c$, $\mathcal{E}_{\beta,K}$ undergoes a continuous bifurcation at $K = K_c^{(2)}(\beta)$ (Theorem 3.1) and that, for $\beta > \beta_c$, $\mathcal{E}_{\beta,K}$ undergoes a discontinuous bifurcation at $K = K_c^{(1)}(\beta)$ (Theorem 3.2). We depict these bifurcations in Figure 8. While the second-order critical values $K_c^{(2)}(\beta)$ are explicitly defined in Theorem 3.6, the first-order critical values $K_c^{(1)}(\beta)$ in the figure are computed numerically. The numerical procedure calculates $K_c^{(1)}(\beta)$ for fixed values of $\beta$ by determining the value of $K$ for which the number of global minimum points of $G_{\beta,K}(z)$ changes from one at $z = 0$ to three at $z = 0$ and $z = \pm\tilde{z}(\beta,K)$, where $\tilde{z}(\beta, K) > 0$. According to these numerical calculations for the discontinuous bifurcation, it appears that $K_c^{(1)}(\beta)$ tends to 1 as $\beta \to \infty$. However, we are unable to prove this limit.

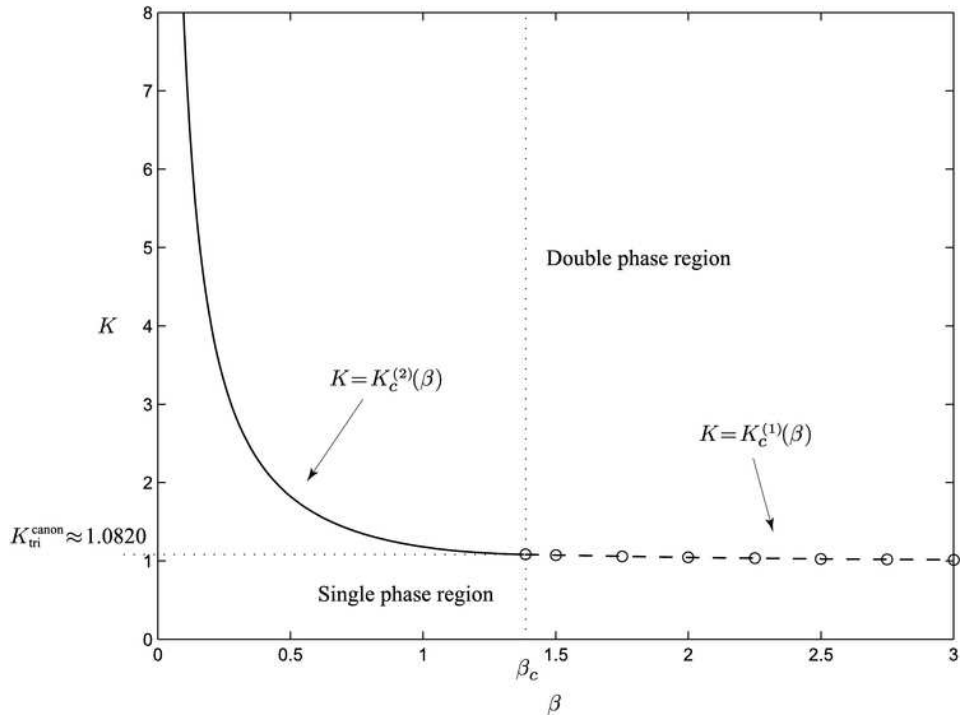

FIG. 8. *Bifurcation diagram for the BEG model with respect to the canonical ensemble.*

In Section 5 we will see that Figure 8 is a phase diagram that describes the phase transitions in the canonical ensemble as $\beta$ changes. We will also show that the nature of the bifurcations studied up to this point by varying $K$, while keeping $\beta$ fixed, is the same if we vary $\beta$ and keep $K$ fixed instead. The latter situation corresponds to what is referred to physically as a phase transition; specifically, the continuous bifurcation corresponds to a second-order phase transition and the discontinuous bifurcation to a first-order phase transition. In order to substantiate this claim concerning the bifurcations and the phase transitions, we have to transfer our analysis of $\mathcal{E}_{\beta,K}$ from fixed $\beta$ and varying $K$ to an analysis of $\mathcal{E}_{\beta,K}$ for fixed $K$ and varying $\beta$.

In the next section we study the BEG model with respect to the microcanonical ensemble.

**4. Structure of the set of microcanonical equilibrium macrostates.** In previous studies of the BEG model with respect to the microcanonical en-semble, results were obtained that either relied on a local analysis or used strictly numerical methods [2, 3, 15]. In this section we provide a global argument to support the existence of a continuous bifurcation exhibited by the set $\mathcal{E}^{u,K}$ of microcanonical equilibrium macrostates for fixed, sufficiently large values of $u$ and for varying $K$. Specifically, for fixed, sufficiently large $u$, $\mathcal{E}^{u,K}$ exhibits a continuous bifurcation as $K$ passes through a critical value $K_c^{(2)}(u)$. The argument is similar to the one employed to analyze the canonical ensemble in Section 3. However, unlike the canonical case, where a rigorous analysis of the structure of the set $\mathcal{E}_{\beta,K}$ of canonical equilibrium macrostates was obtained for all values of $\beta$ and $K$, the analysis of $\mathcal{E}^{u,K}$ for sufficiently large $u$ and varying $K$ relies on a mixture of analysis and numerical methods. At the end of this section we summarize the numerical methods used to deduce the existence of a discontinuous bifurcation exhibited by $\mathcal{E}^{u,K}$ for fixed, sufficiently small $u$ and varying $K$. In Section 5 we show how to extrapolate this information to information concerning the phase transition behavior of the microcanonical ensemble for varying $u$: a continuous, second-order phase transition for all sufficiently large values of $K$ and a discontinuous, first-order phase transition for all sufficiently small values of $K$.

We begin by recalling several definitions from Section 2. $\mathcal{P}(\Lambda)$ denotes the set of probability measures with support $\Lambda = \{-1, 0, 1\}$; $\rho$ denotes the measure $\frac{1}{3}(\delta_{-1} + \delta_0 + \delta_1) \in \mathcal{P}(\Lambda)$; for $\mu \in \mathcal{P}(\Lambda)$,

$$R(\mu|\rho) = \sum_{i=-1}^{1} \mu_i \log 3\mu_i$$

denotes the relative entropy of $\mu$ with respect to $\rho$; and $f_K(\mu)$ is defined by

$$f_K(\mu) = \int_\Lambda y^2 \mu(dy) - K\left(\int_\Lambda y\mu(dy)\right)^2$$
$$= (\mu_1 + \mu_{-1}) - K(\mu_1 - \mu_{-1})^2.$$

For $K > 0$, we also defined the set of microcanonical equilibrium macrostates by

$$\begin{aligned}\mathcal{E}^{u,K} &= \{\nu \in \mathcal{P}(\Lambda) : I^{u,K}(\nu) = 0\} \\ &= \{\nu \in \mathcal{P}(\Lambda) : \nu \text{ minimizes } R(\nu|\rho) \text{ subject to } f_K(\nu) = u\},\end{aligned} \quad (4.1)$$

$\mathcal{E}^{u,K}$ is well defined for $K > 0$ and $u \in \operatorname{dom} s_K = [\min(1-K, 0), 1]$. Throughout this section we fix $u \in \operatorname{dom} s_K$; $s_K$ is defined in (2.4).

Determining the elements in $\mathcal{E}^{u,K}$ requires solving a constrained minimization problem, which is the dual of the unconstrained minimization problem associated with the set $\mathcal{E}_{\beta,K}$ of canonical equilibrium macrostates defined in (2.7). In order to simplify the analysis of the set $\mathcal{E}^{u,K}$, we employ the technique used in [2] to reduce the constrained minimization problem defining $\mathcal{E}^{u,K}$ to another minimization problem that is more easily studied. For fixed $K > 0$ and $u \in \operatorname{dom} s_K$, we define

$$\mathcal{D}_{u,K} = \{\mu \in \mathcal{P}(\Lambda) : f_K(\mu) = u\}. \tag{4.2}$$

For $\mu \in \mathcal{D}_{u,K}$, let $z = \mu_1 - \mu_{-1}$ and $q = \mu_1 + \mu_{-1}$. Since $\mu \in \mathcal{D}_{u,K}$ implies that

$$f_K(\mu) = (\mu_1 + \mu_{-1}) - K(\mu_1 - \mu_{-1})^2 = u,$$

we see that $q = u + Kz^2$. Thus, for $\mu \in \mathcal{D}_{u,K}$, we have

$$\begin{aligned}
R(\mu|\rho) &= \sum_{i=-1}^{1} \mu_i \log 3\mu_i \\
&= \frac{q-z}{2} \log\left[\frac{3}{2}(q-z)\right] + \frac{q+z}{2} \log\left[\frac{3}{2}(q+z)\right] \\
&\quad + (1-q)\log[3(1-q)] \\
&= \frac{q+z}{2} \log(q+z) + \frac{q-z}{2} \log(q-z) \\
&\quad + (1-q)\log(1-q) - (q\log 2 - \log 3).
\end{aligned}$$

Setting $q = u + Kz^2$, we define the quantity

$$\begin{aligned}
R_{u,K}(z) &= \frac{q+z}{2} \log(q+z) + \frac{q-z}{2} \log(q-z) \\
&\quad + (1-q)\log(1-q) - (q\log 2 - \log 3)
\end{aligned} \tag{4.3}$$

and the set

$$\mathcal{M}_{u,K} = \{z \in \mathbb{R} : z = \mu_1 - \mu_{-1} \text{ for some } \mu \in \mathcal{D}_{u,K}\}. \tag{4.4}$$

The derivation of $R_{u,K}$ makes it clear that $\mathcal{M}_{u,K} \subset (-1,1)$ is the domain of $R_{u,K}$.

We next introduce the set

$$\tilde{\mathcal{E}}^{u,K} = \{\tilde{z} \in \mathcal{M}_{u,K} : \tilde{z} \text{ minimizes } R_{u,K}(z)\}.$$

The following theorem states a one-to-one correspondence between the elements of $\mathcal{E}^{u,K}$ and $\tilde{\mathcal{E}}^{u,K}$ under an assumption on the structure of $\tilde{\mathcal{E}}^{u,K}$. In [15], for particular values of $u$ and $K$, numerical experiments show that $\tilde{\mathcal{E}}^{u,K}$ consists of either $1, 2$ or $3$ points. Although we are not able to prove that this is valid for all $u \in \operatorname{dom} s_K$ and $K > 0$, because of our numerical computations, we make it an assumption in the next theorem.

PHASE TRANSITIONS IN THE MEAN-FIELD BEG MODEL 39

THEOREM 4.1. *Fix $K > 0$ and $u \in \mathrm{dom}\, s_K$. Suppose that $\tilde{\mathcal{E}}^{u,K} = \{z_\alpha\}_{\alpha=1}^r$, where $r$ equals 1, 2 or 3. Define $\nu_\alpha = \sum_{i=-1}^1 \nu_{\alpha,i}\delta_i \in \mathcal{P}(\Lambda)$ by the formulas*

$$\nu_{\alpha,1} = \frac{u + Kz_\alpha^2 + z_\alpha}{2}, \qquad \nu_{\alpha,-1} = \frac{u + Kz_\alpha^2 - z_\alpha}{2}, \qquad \nu_{\alpha,0} = 1 - \nu_{\alpha,1} - \nu_{\alpha,-1}.$$

*Then $\mathcal{E}^{u,K}$ consists of the distinct elements $\nu_\alpha, \alpha = 1, \ldots, r$.*

PROOF. Using the definition (4.2) of $\mathcal{D}_{u,K}$, we can rewrite the set $\mathcal{E}^{u,K}$ of microcanonical equilibrium macrostates defined in (4.1) as

$$\mathcal{E}^{u,K} = \{\nu \in \mathcal{D}_{u,K} : \nu \text{ is a minimum point of } R(\mu|\rho)\}.$$

We show that, for $\alpha = 1, \ldots, r$, $f_K(\nu_\alpha) = u$ and $R(\nu_\alpha|\rho) < R(\mu|\rho)$ for all $\mu \in \mathcal{D}_{u,K}$ for which $\mu \neq \nu_\alpha$.

From the definition of $\nu_\alpha$, we have

$$f_K(\nu_\alpha) = (\nu_{\alpha,1} + \nu_{\alpha,-1}) - K(\nu_{\alpha,1} - \nu_{\alpha,-1})^2 = (u + Kz_\alpha^2) - Kz_\alpha^2 = u.$$

Therefore, $\nu_\alpha \in \mathcal{D}_{u,K}$ for all $\alpha = 1, \ldots, r$. Since for all $z_\alpha, z_\ell \in \tilde{\mathcal{E}}^{u,K}$, $\alpha, \ell = 1, \ldots, r$,

$$R(\nu_\alpha|\rho) = R_{u,K}(z_\alpha) = R_{u,K}(z_\ell) = R(\nu_\ell|\rho),$$

it follows that $R(\nu_\alpha|\rho)$ are equal for all $\alpha = 1, \ldots, r$.

We now consider $\mu = \sum_{i=-1}^1 \mu_i \delta_i \in \mathcal{D}_{u,K}$ such that $\mu \neq \nu_\alpha$ for all $\alpha = 1, \ldots, r$. Defining $\zeta = \mu_1 - \mu_{-1}$, we claim that $\zeta \neq z_\alpha$ for all $\alpha = 1, \ldots, r$. Suppose otherwise; that is, for some $z_\alpha$,

(4.5) $$\mu_1 - \mu_{-1} = \zeta = z_\alpha = \nu_{\alpha,1} - \nu_{\alpha,-1}.$$

But $\mu \in \mathcal{D}_{u,K}$ implies that $f_K(\mu) = u = f_K(\nu_\alpha)$ and, thus, that

(4.6) $$\mu_1 + \mu_{-1} = \nu_{\alpha,1} + \nu_{\alpha,-1}.$$

Combining (4.5) and (4.6) yields the contradiction that $\mu = \nu_\alpha$. Because $\zeta \neq z_\alpha$ for all $\alpha = 1, \ldots, r$, it follows that $\zeta \notin \tilde{\mathcal{E}}^{u,K}$ and, thus, that $R_{u,K}(z_\alpha) < R_{u,K}(\zeta)$ for all $\alpha = 1, \ldots, r$. As a result, for $\alpha = 1, \ldots, r$, we have

$$R(\nu_\alpha|\rho) = R_{u,K}(z_\alpha) < R_{u,K}(\zeta) = R(\mu|\rho).$$

We complete the proof by showing that if $z_\alpha \neq z_\ell$, then $\nu_\alpha \neq \nu_\ell$. Indeed, if $\nu_\alpha = \nu_\ell$, then, for each choice of sign, we would have $Kz_\alpha^2 \pm z_\alpha = Kz_\ell^2 \pm z_\ell$. Since this leads to the contradiction that $z_\alpha = z_\ell$, the proof of the theorem is complete.

□

Theorem 4.1 allows us to analyze the set $\mathcal{E}^{u,K}$ of microcanonical equilibrium macrostates by calculating the minimum points of the function $R_{u,K}$ defined in (4.3). Define

$$\varphi_{u,K}(z) = \frac{q+z}{2}\log(q+z) + \frac{q-z}{2}\log(q-z) + (1-q)\log(1-q),$$

where $q = u + Kz^2$. With this notation (4.3) becomes

$$R_{u,K}(z) = \varphi_{u,K}(z) - (u + Kz^2)\log 2 + \log 3.$$

This separation of $R_{u,K}$ into the nonlinear component $\varphi_{u,K}$ and the quadratic component is similar to the method used in Sections 3.2 and 3.3 in determining the elements in the set $\tilde{\mathcal{E}}_{\beta,K}$. There we separated the minimizing function $F_{\beta,K}(w)$ into a nonlinear component $c_\beta(w)$ and a quadratic component $w^2/(4\beta K)$; minimum points of $F_{\beta,K}$ satisfy $F'_{\beta,K}(w) = c'_\beta(w) - w/(2\beta K) = 0$. Solving this equation was greatly facilitated by understanding the concavity and convexity properties of $c_\beta$, which are proved in Theorem 3.5.

Following the success of this method in studying the canonical ensemble, we apply a similar technique to determine the minimum points of $R_{u,K}$. We call a pair $(u, K)$ admissible if $u \in \mathrm{dom}\, s_K$. While an analytic proof could not be found, our numerical experiments show that there exists a curve $K = C(u)$ in the $(u, K)$-plane such that for all admissible $(u, K)$ lying above the graph of this curve, $\varphi'_{u,K}$ is strictly convex on its positive domain. The graph of $K = C(u)$ is depicted in Figure 9. We denote by $G^+$ the set of admissible $(u, K)$ lying above this graph and by $G_-$ the set of admissible $(u, K)$ lying below this graph. Using a similar argument as in the proof of Theorem 3.6 for the canonical case, we are led to believe that, for all $(u, K) \in G^+$, the BEG model with respect to the microcanonical ensemble exhibits a continuous bifurcation in $K$; that is, there exists a critical value $K_c^{(2)}(u) > 0$ such that the following hold:

- For $0 < K \leq K_c^{(2)}(u)$, $\tilde{\mathcal{E}}^{u,K} = \{0\}$.
- For $K > K_c^{(2)}(u)$, there exists a positive number $\tilde{z}(u, K)$ such that $\tilde{\mathcal{E}}^{u,K} = \{\pm \tilde{z}(u, K)\}$.
- $\lim_{K \to (K_c^{(2)}(u))^+} \tilde{z}(u, K) = 0$.

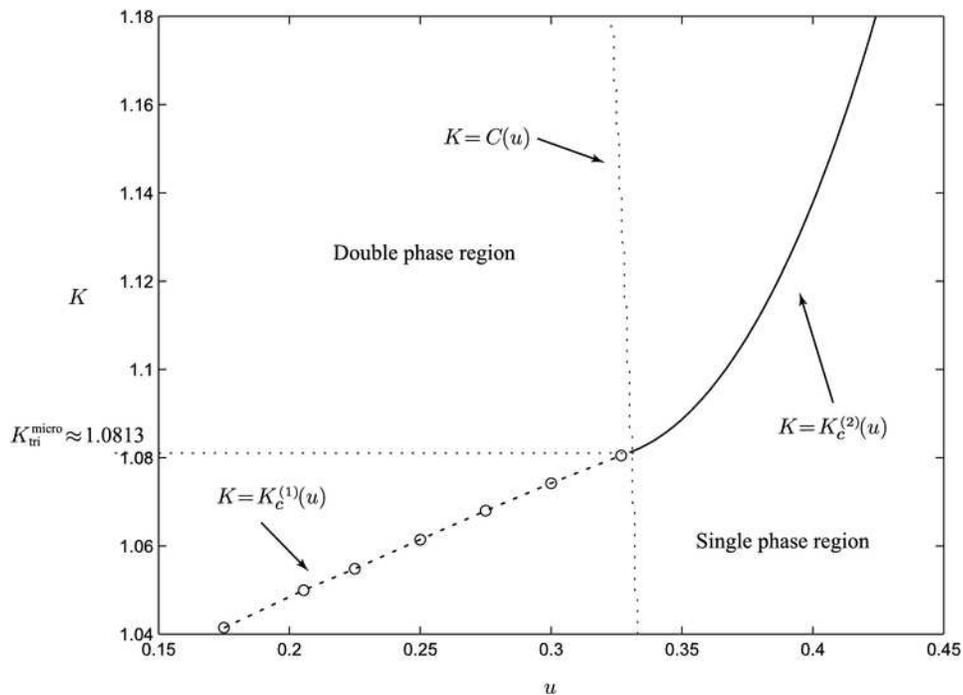

FIG. 9.   *Bifurcation diagram for the BEG model with respect to the microcanonical ensemble.*

Combined with the one-to-one correspondence between the elements of $\tilde{\mathcal{E}}^{u,K}$ and $\mathcal{E}^{u,K}$ proved in Theorem 4.1, the structure of $\tilde{\mathcal{E}}^{u,K}$ just given yields a continuous bifurcation in $K$ exhibited by $\mathcal{E}^{u,K}$ for $(u, K)$ lying in the region $G^+$ above the graph of the curve $K = C(u)$. Similar to the definition of the critical value $K_c^{(2)}(\beta)$ given in (3.19) for the continuous bifurcation in $K$ exhibited by $\tilde{\mathcal{E}}_{\beta,K}$, the critical value $K_c^{(2)}(u)$ is the solution of the equation

$$R''_{u,K}(0) = 0 \quad \text{or} \quad \varphi''_{u,K}(0) = 2K \log 2.$$

Consequently, since $\varphi''_{u,K}(0) = 1/u + 2K[\log(u/(1-u))]$, we define the second-order critical value to be

$$(4.7) \qquad K_c^{(2)}(u) = \frac{\varphi''_{u,K}(0)}{2 \log 2} = \frac{1}{2u \log(2(1-u)/u)}.$$

The derivation of this formula for $K_c^{(2)}(u)$ for the critical values of the continuous bifurcation in $K$ exhibited by $\mathcal{E}^{u,K}$ rests on the existence of the curve $K = C(u)$, which in turn was derived numerically. However, the accuracy of (4.7) is supported by the fact that the graph of the curve $K_c^{(2)}(u)$ fits the critical values derived numerically in Figures 2 and 3 of [15].

For values of $(u, K)$ lying in the region $G_-$ below the graph of the curve $K = C(u)$, the strict convexity behavior of $\varphi'_{u,K}$ no longer holds. Therefore, numerical computations were used to determine the behavior of $R_{u,K}$ for such $(u, K)$, showing a discontinuous bifurcation in $K$ in this region. Specifically, there exists a critical value $K_c^{(1)}(u)$ such that the following hold:

- For $0 < K < K_c^{(1)}(u)$, $\tilde{\mathcal{E}}^{u,K} = \{0\}$.
- For $K = K_c^{(1)}(u)$, there exists $\tilde{z}(u, K) > 0$ such that $\tilde{\mathcal{E}}^{u,K} = \{0, \pm \tilde{z}(u, K)\}$.
- For $K > K_c^{(1)}(u)$, there exists $\tilde{z}(u, K) > 0$ such that $\tilde{\mathcal{E}}^{u,K} = \{\pm \tilde{z}(u, K)\}$.

The critical values $K_c^{(1)}(u)$ were computed numerically by determining the value of $K$ for which the number of global minimum points of $R_{u,K}(z)$ changes from one at $z = 0$ to three at $z = 0$ and $z = \pm \tilde{z}(u, K)$, $\tilde{z}(u, K) > 0$.

The results of this section are summarized in the bifurcation diagram for the BEG model with respect to the microcanonical ensemble, which appears in Figure 9. In the next section we will see that Figure 9 is a phase diagram that describes the phase transition in the microcanonical ensemble as $u$ changes. In order to substantiate this, we have to transfer our analysis of $\mathcal{E}^{u,K}$ from fixed $u$ and varying $K$ to an analysis of $\mathcal{E}^{u,K}$ for fixed $K$ and varying $u$.

**5. Comparison of phase diagrams for the two ensembles.** We end our analysis of the canonical and microcanonical ensembles by explaining what our results imply concerning the nature of the phase transitions in the BEG model. These phase transitions are defined by varying $\beta$ and $u$, the two parameters that define the ensembles. As we will see, the order of the phase transitions is a structural property of the phase diagram in the sense that it is the same whether we vary $K$ or $\beta$ in the canonical ensemble and $K$ or $u$ in the microcanonical ensemble while keeping the other parameter fixed.

Before doing this, we first review one of the main contributions of the preceding two sections, which is to analyze the bifurcation behavior of the sets $\mathcal{E}_{\beta,K}$ and $\mathcal{E}^{u,K}$ of equilibrium macrostates with respect to both the canonical and microcanonical ensembles. Figure 8 summarizes the canonical analysis and Figure 9 the microcanonical analysis. The figures exhibit two different values of $K$ called tricritical values and denoted by $K_{\text{tri}}^{\text{canon}}$ and $K_{\text{tri}}^{\text{micro}}$. As we soon explain, at each of these values of $K$ the corresponding ensemble changes its behavior from a continuous, second-order phase transition to a discontinuous, first-order phase transition.

For the canonical ensemble, the tricritical value in Figure 8 is given by

$$K_{\text{tri}}^{\text{canon}} = K_c^{(2)}(\beta_c) = K_c^{(2)}(\log 4) \approx 1.0820,$$

where $K_c^{(2)}(\beta)$ is defined in (3.19). With respect to the microcanonical ensemble, the tricritical value $K_{\text{tri}}^{\text{micro}}$ is the value of $K$ at which the curves $K = C(u)$ and $K_c^{(2)}(u)$ shown in Figure 9 intersect. From the numerical calculation of the curve $K = C(u)$, we obtain the following approximation for the tricritical value $K_{\text{tri}}^{\text{micro}}$:

$$K_{\text{tri}}^{\text{micro}} \approx 1.0813.$$

These values of $K_{\text{tri}}^{\text{canon}}$ and $K_{\text{tri}}^{\text{micro}}$ agree with the values derived in [2] via a local analysis and numerical computations.

We first illustrate how our analysis of $\mathcal{E}_{\beta,K}$ in Theorems 3.1 and 3.2 for fixed $\beta$ and varying $K$ yields a continuous, second-order phase transition and a discontinuous, first-order phase transition with respect to the canonical ensemble. These phase transitions are defined for fixed $K$ and varying $\beta$, the thermodynamic parameter that defines the ensemble. In order to study the phase transition, we must therefore transform the analysis of $\mathcal{E}_{\beta,K}$ for fixed $\beta$ and varying $K$ to an analysis of the same set for fixed $K$ and varying $\beta$. After we consider the microcanonical phase transition in an analogous way, we will focus on the region

$$K_{\text{tri}}^{\text{micro}} \approx 1.0813 < K < 1.0820 \approx K_{\text{tri}}^{\text{canon}}.$$

As we will point out, the fact that for $K$ in this region the two ensembles exhibit different phase transition behavior—discontinuous for the canonical and continuous for the microcanonical—is closely related to the phenomenon of ensemble nonequivalence in the model.

We begin with the continuous phase transition for the canonical ensemble. Figure 8 exhibits a monotonically decreasing function $K = K_c^{(2)}(\beta)$ for $0 < \beta < \beta_c = \log 4$. Inverting this function yields a monotonically decreasing function $\beta = \beta_c^{(2)}(K)$ for $K > K_{\text{tri}}^{\text{canon}} = K_c^{(2)}(\beta_c) \approx 1.0820$. Consider, for fixed $K > K_{\text{tri}}^{\text{canon}}$ and small $\delta > 0$, values of $\beta \in (\beta_c^{(2)}(K) - \delta, \beta_c^{(2)}(K) + \delta)$. Our analysis of $\mathcal{E}_{\beta,K}$ in Theorem 3.1 shows the following:

- For $\beta \in (\beta_c^{(2)}(K) - \delta, \beta_c^{(2)}(K)]$, the model exhibits a single phase $\rho_\beta$.
- For $\beta \in (\beta_c^{(2)}(K), \beta_c^{(2)}(K) + \delta)$, the model exhibits two distinct phases $\nu^+(\beta, K)$ and $\nu^-(\beta, K)$.

We claim that, for fixed $K > K_{\text{tri}}^{\text{canon}}$, this is a second-order phase transition; that is, as $\beta \to (\beta_c^{(2)}(K))^+$, we have $\nu^+(\beta, K) \to \rho_\beta$ and $\nu^-(\beta, K) \to \rho_\beta$. To see this, we recall from Figure 1(b) that, for $\beta = \beta_c^{(2)}(K)$, the graph of the linear component $w/(2\beta K)$ of $F'_{\beta,K}(w)$ is tangent to the graph of the nonlinear component $c'_\beta(w)$ of $F'_{\beta,K}(w)$ at the origin. This figure was referred to in Section 3.1 when we analyzed the structure of the set $\tilde{\mathcal{E}}_{\beta,K}$ (Theorem 3.6). Since both components of $F'_{\beta,K}(w)$ are continuous with respect to $\beta$, a perturbation in $\beta$ yields a continuous phase transition in $\tilde{\mathcal{E}}_{\beta,K}$ and thus in $\mathcal{E}_{\beta,K}$. A similar argument shows that each of the double phases $\nu^+(\beta, K)$ and $\nu^-(\beta, K)$ are continuous functions of $\beta$ for $\beta > \beta_c^{(2)}(K)$.

We now analyze the discontinuous phase transition for the canonical ensemble in a similar way. Figure 8 exhibits a monotonically decreasing function $K = K_c^{(1)}(\beta)$ for $\beta > \beta_c = \log 4$. Inverting this function yields a monotonically decreasing function $\beta = \beta_c^{(1)}(K)$ for $0 < K < K_{\text{tri}}^{\text{canon}} \approx 1.0820$. For fixed $0 < K < K_{\text{tri}}^{\text{canon}}$ and small $\delta > 0$, consider values of $\beta \in (\beta_c^{(1)}(K) - \delta, \beta_c^{(1)}(K) + \delta)$. Our analysis of $\mathcal{E}_{\beta,K}$ in Theorem 3.2 shows the following:

- For $\beta \in (\beta_c^{(1)}(K) - \delta, \beta_c^{(1)}(K))$, the model exhibits a single phase $\rho_\beta$.
- For $\beta = \beta_c^{(1)}(K)$, the model exhibits three distinct phases $\rho_\beta$, $\nu^+(\beta, K)$, and $\nu^-(\beta, K)$.

- For $\beta \in (\beta_c^{(1)}(K), \beta_c^{(1)}(K) + \delta)$, the model exhibits two distinct phases $\nu^+(\beta, K)$ and $\nu^-(\beta, K)$.

We claim that, for fixed $0 < K < K_{\text{tri}}^{\text{canon}}$, this is a first-order phase transition; that is, as $\beta \to (\beta_c^{(1)}(K))^+$, we have, for each choice of sign, $\nu^\pm(\beta, K) \to \nu^\pm(\beta_c^{(1)}(K), K) \neq \rho_\beta$. To see this, we recall from Figure 7(a) that, for $\beta = \beta_c^{(1)}(K)$, the graph of the linear component $w/(2\beta K)$ of $F'_{\beta,K}(w)$ intersects the graph of the nonlinear component $c'_\beta(w)$ of $F'_{\beta,K}(w)$ in five places such that the signed area between the two graphs is 0. This results in three values of $w$ that are global minimum points of $F_{\beta,K}$; namely, $w = 0, \tilde{w}(\beta, K), -\tilde{w}(\beta, K)$ (Theorem 3.8). These three values of $w$ give rise to three values of $z = w/(2\beta K)$ that constitute the set $\tilde{\mathcal{E}}_{\beta,K}$ for $\beta = \beta_c^{(2)}(K)$. Since both components of $F'_{\beta,K}(w)$ are continuous with respect to $\beta$, a perturbation in $\beta$ yields a discontinuous phase transition in $\tilde{\mathcal{E}}_{\beta,K}$ and thus in $\mathcal{E}_{\beta,K}$. A similar argument shows that each of the equilibrium macrostates $\nu^+(\beta, K)$ and $\nu^-(\beta, K)$ are continuous functions of $\beta$ for $\beta > \beta_c^{(2)}(K)$.

The phase transitions for the microcanonical ensemble are defined for fixed $K$ and varying $u$, the thermodynamic parameter defining the ensemble. Therefore, in order to study these phase transitions, we must transform the analysis of $\mathcal{E}^{u,K}$ done in Section 4 for fixed $u$ and varying $K$ to an analysis of the same set for fixed $K$ and varying $u$. This is carried out in a way that is similar to what we have just done for the canonical ensemble. In particular, we find that, for $K > K_{\text{tri}}^{\text{micro}} \approx 1.0813$, the BEG model with respect to the microcanonical ensemble exhibits a continuous, second-order phase transition and that, for $0 < K < K_{\text{tri}}^{\text{micro}}$, the model exhibits a discontinuous, first-order phase transition.

We now focus on values of $K$ satisfying $K_{\text{tri}}^{\text{micro}} < K < K_{\text{tri}}^{\text{canon}}$. As we have just seen, for such $K$, the two ensembles exhibit different phase transition behavior: for $K_{\text{tri}}^{\text{micro}} < K$, the microcanonical ensemble undergoes a continuous, second-order phase transition, while for $0 < K < K_{\text{tri}}^{\text{canon}}$, the canonical ensemble undergoes a discontinuous, first-order phase transition. This observation is consistent with a numerical calculation given in Figure 10 showing that, for a fixed value of $K \in (K_{\text{tri}}^{\text{micro}}, K_{\text{tri}}^{\text{canon}})$, there exists a subset of the microcanonical equilibrium macrostates that are not realized canonically [15]. As a result, for this value of $K$, the two ensembles are nonequivalent at the level of equilibrium macrostates.

Figures 10(a) and 10(b) exhibit, for a range of values of $u$ and $\beta$, the structure of the set $\mathcal{E}^{u,K}$ of microcanonical equilibrium macrostates and the set $\mathcal{E}_{\beta,K}$ of canonical equilibrium macrostates for $K = 1.0817$. This value of $K$ lies in the interval $(K_{\text{tri}}^{\text{micro}}, K_{\text{tri}}^{\text{canon}}) \approx (1.0813, 1.0820)$. Each equilibrium macrostate in $\mathcal{E}^{u,K}$ and $\mathcal{E}_{\beta,K}$ is an empirical measure having the form

$$\nu = \nu_1 \delta_1 + \nu_0 \delta_0 + \nu_{-1} \delta_{-1}.$$

In both figures the solid and dashed curves can be taken to represent the components $\nu_1$ and $\nu_{-1}$. The components $\nu_1$ and $\nu_{-1}$ in the microcanonical ensemble are functions of $u$ [Figure 10(a)] and in the canonical ensemble are functions of $\beta$ [Figure 10(b)]. Figures 10(a) and 10(b) were taken from [15].

Comparing the two figures reveals that the ensembles are nonequivalent for this value of $K$. Specifically, because of the discontinuous, first-order phase transition in the canonical ensemble, there exists a subset of $\mathcal{P}(\Lambda)$ that is not realized by $\mathcal{E}_{\beta,K}$ for any $\beta > 0$. On the other hand, since the set $\mathcal{E}^{u,K}$ of microcanonical equilibrium macrostates exhibits a continuous, second-order phase transition, the subset of $\mathcal{P}(\Lambda)$ not realized canonically is realized microcanonically. As a result, there exists a nonequivalence of ensembles at the level of equilibrium macrostates. The reader is referred to [15] for a more complete analysis of ensemble equivalence and nonequivalence for the BEG model.

**6. Limit theorems for the total spin with respect to $P_{n,\beta,K}$.** In Section 3.1 we rewrote the canonical ensemble $P_{n,\beta,K}$ for the BEG model in terms of the total spin $S_n$. This allowed us to reduce the analysis of the set $\mathcal{E}_{\beta,K}$ of canonical equilibrium macrostates to that of a Curie–Weiss-type model. We end this paper by deriving limit theorems for the $P_{n,\beta,K}$-distributions of appropriately scaled partial sums $S_n = \sum_{j=1}^n \omega_j$, which represents the total spin in the model. Since $S_n/n = \int_{\{-1,0,1\}} y L_n(dy)$, the limit theorems for $S_n$ are also limit theorems for the empirical measures $L_n$. As we will see, the new limit theorems follow from those for the Curie–Weiss model proved in [12, 14].

Let $\tau$ be a Borel probability measure on $\mathbb{R}$ satisfying $\int_{\mathbb{R}} \exp[bx^2] \tau(dx) < \infty$ for all $b > 0$. The Curie–Weiss model considered in [12, 14] is defined in terms

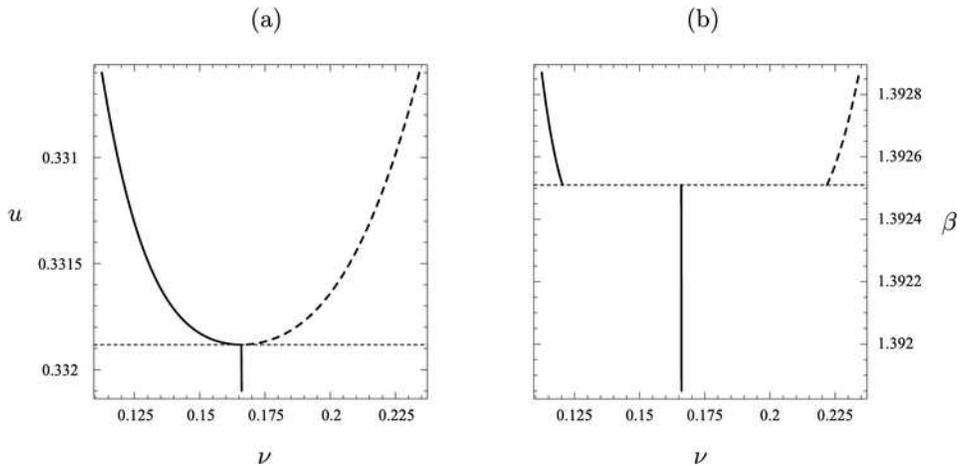

FIG. 10. *Structure of* (a) *the set* $\mathcal{E}^{u,K}$ *and* (b) *the set* $\mathcal{E}_{\beta,K}$ *for* $K = 1.0817$.

of a canonical ensemble on $(\mathbb{R}^n, \mathcal{B}_{\mathbb{R}_n})$ given by

$$(6.1) \qquad P_{n,\beta}^\tau(d\omega) = \frac{1}{Z_n^\tau(\beta)} \cdot \exp\left[\frac{n\beta}{2}\left(\frac{S_n(\omega)}{n}\right)^2\right] P_n^\tau(d\omega).$$

In this formula $\beta > 0$, $P_n^\tau$ is the product measure on $\mathbb{R}^n$ with identical one-dimensional marginals $\tau$, and $Z_n^\tau(\beta)$ is a normalization making $P_{n,\beta}^\tau$ a probability measure. The canonical ensemble for the BEG model is defined by the measure $P_{n,\beta,K}$ in (2.1), which is re-expressed in (3.2) as a Curie–Weiss-type measure. This measure has the form (6.1), in which $\beta$ is replaced by $2\beta K$ and $\tau$ equals the measure $\rho_\beta$ defined in (3.1).

For $t \in \mathbb{R}$, define $c^\tau(t) = \log \int_\mathbb{R} \exp(t\omega_1)\tau(d\omega_1)$. As shown in [12, 14], the $P_{n,\beta}^\tau$-limits for $S_n$ are determined by the global minimum points of the function

$$(6.2) \qquad G_\beta^\tau(z) = \tfrac{1}{2}\beta z^2 - c^\tau(\beta z).$$

Let $\tilde{z}$ be a global minimum point of $G_\beta^\tau$. Since $G_\beta^\tau$ is real analytic, there exists a positive integer $r = r(\tilde{z})$ such that $(G_\beta^\tau)^{(2r)}(\tilde{z}) > 0$ and

$$G_\beta^\tau(z) = G_\beta^\tau(\tilde{z}) + \frac{(G_\beta^\tau)^{(2r)}(\tilde{z})}{(2r)!}(z - \tilde{z})^{2r} + O((z - \tilde{z})^{2r+1}) \qquad \text{as } z \longrightarrow \tilde{z}.$$

We call $r(\tilde{z})$ the type of the minimum point $\tilde{z}$. If $r = 1$, then $(G_\beta^\tau)''(\tilde{z}) = \beta - \beta^2 (c_\beta^\tau)''(\tilde{z})$, and if $r \geq 2$, then $(G_\beta^\tau)^{(2r)}(\tilde{z}) = -\beta^{2r}(c_\beta^\tau)^{(2r)}(\tilde{z})$.

The canonical ensemble $P_{n,\beta,K}$ for the BEG model has the form of the Curie–Weiss measure $P_{n,\beta}^\tau$ with $\beta$ replaced by $2\beta K$ and $\tau = \rho_\beta$. Therefore, the function that plays the role of $G_\beta^\tau$ for the BEG model is $G_{2\beta K}^{\rho_\beta}$. This coincides with the function

$$G_{\beta,K}(z) = \beta K z^2 - c_\beta(2\beta K z)$$
$$= \beta K z^2 - \log \int_{\{-1,0,1\}} \exp(2\beta K \omega_1) \rho_\beta(d\omega_1),$$

defined in (3.8). For $0 < \beta \leq \beta_c$ and $K > 0$, detailed information about the set $\tilde{\mathcal{E}}_{\beta,K}$ of global minimum points of $G_{\beta,K}$ is given in Theorem 3.6; for $\beta > \beta_c$ and $K > 0$, detailed information about $\tilde{\mathcal{E}}_{\beta,K}$ is given in Theorem 3.8.

We next indicate the form of the limit theorems for the Curie–Weiss model, restricting to those cases that arise in the BEG model. The first, Theorem 6.1, states limits that are valid when $G_\beta^\tau$ has a unique global minimum point at $z = 0$. The second, Theorem 6.2, states a conditioned limit that is valid when $G_\beta^\tau$ has multiple global minimum points all of type 1.

A law of large numbers for $S_n/n$ is given in part (a) of Theorem 6.1. In part (b) $f_{0,\sigma^2(\beta)}$ denotes the density of a $N(0, \sigma^2(\beta))$ random variable with

$$(6.3) \qquad \sigma^2(\beta) = \frac{\beta \cdot (c_\beta^\tau)''(0)}{(G_\beta^\tau)''(0)}.$$

When the type of the minimum point at 0 is $r = 1$, $\sigma^2(\beta) > 0$ because, in this case, $(G_\beta^\tau)''(0) > 0$ and, in general, $(c_\beta^\tau)''(0) > 0$. If $f$ is a nonnegative, integrable function on $\mathbb{R}$, then, for $r \in \mathbb{N}$, we write

$$P_{n,\beta}^\tau\{S_n/n^{1-1/2r} \in dx\} \implies f(x)\,dx$$

to mean that, as $n \to \infty$, the $P_{n,\beta}^\tau$-distributions of $S_n/n^{1-1/2r}$ converge weakly to a distribution having a density proportional to $f$. When $r = 1$, $f = f_{0,\sigma^2(\beta)}$, and the limit is a central-limit-type theorem with scaling $n^{1/2}$. When $r \geq 2$, the limits involve the nonclassical scaling $n^{1-1/2r}$, and the $P_{n,\beta}^\tau$-distributions of the scaled random variables converge weakly to a distribution having a density proportional to $\exp[-\text{const} \cdot x^{2r}]$. Theorem 6.1 is proved in Theorem 2.1 in [12] for $\beta = 1$; rescaling yields the more general form given here.

THEOREM 6.1. *Consider the Curie–Weiss model, for which the canonical ensemble $P_{n,\beta}^\tau$ is defined by (6.1). For $\beta > 0$, assume that $G_\beta^\tau$ has a unique global minimum point at $z = 0$ having type $r$. Let $f_{0,\sigma^2(\beta)}$ be the density of a $N(0, \sigma^2(\beta))$ random variable, where $\sigma^2(\beta)$ is the positive quantity defined in (6.3). The following conclusions hold:*

(a) $P_{n,\beta}^\tau\{S_n/n \in dx\} \Longrightarrow \delta_0$ *as* $n \to \infty$.
(b) *As* $n \to \infty$,

$$P_{n,\beta}^\tau\left\{\frac{S_n}{n^{1-1/2r}} \in dx\right\} \implies \begin{cases} f_{0,\sigma^2(\beta)}(x)\,dx, & \text{for } r = 1, \\ \exp(-(G_\beta^\tau)^{(2r)}(0) \cdot x^{2r}/(2r)!)\,dx, & \text{for } r \geq 2. \end{cases}$$

The next theorem is valid when $G_\beta^\tau$ has multiple global minimum points all of type 1. Part (a), proved in Theorem 3.8 in [12], states a law of large numbers for $S_n/n$. Part (b), proved in Theorem 2.4 in [14], states a conditioned limit. For each global minimum point $\tilde{z}$ of type 1, we define the positive quantity

(6.4) $$\sigma^2(\beta, \tilde{z}) = \frac{\beta \cdot (c_\beta^\tau)''(\beta \tilde{z})}{(G_\beta^\tau)''(\beta \tilde{z})}.$$

THEOREM 6.2. *Consider the Curie–Weiss model, for which the canonical ensemble $P_{n,\beta}^\tau$ is defined by (6.1). For $\beta > 0$, assume that $G_\beta^\tau$ has global minimum points, all of type 1, at $\{z_1, \ldots, z_m\}$ for $m \geq 2$. For each $j = 1, \ldots, m$, we define*

$$b_j = \frac{\sigma^2(\beta, z_j)}{\sum_{\ell=1}^m \sigma^2(\beta, z_\ell)},$$

where $\sigma^2(\beta, z_j)$ is the positive quantity defined in (6.4). Let $f_{0,\sigma^2(\beta,z_j)}$ be the density of a $N(0, \sigma^2(\beta, z_j))$ random variable. The following conclusions hold:

(a) $P_{n,\beta}^\tau \{S_n/n \in dx\} \Longrightarrow \sum_{j=1}^m b_j \delta_{z_j}$ as $n \to \infty$.
(b) There exists $\alpha = \alpha(z_j) > 0$ such that, for any $a \in (0, \alpha)$,

$$P_{n,\beta}^\tau \left\{ \frac{S_n - nz_j}{n^{1/2}} \in dx \bigg| \frac{S_n}{n} \in [z_j - a, z_j + a] \right\}$$
$$\Longrightarrow f_{0,\sigma^2(\beta,z_j)}(x)\, dx \quad \text{as } n \to \infty.$$

In order to adapt these limit theorems to the BEG model, we now classify each of the points in $\tilde{\mathcal{E}}_{\beta,K}$ by type. $\tilde{\mathcal{E}}_{\beta,K}$ denotes the set of global minimum points of $G_{\beta,K} = G_{2\beta K}^{\rho_\beta}$, which plays the same role for the BEG model as $G_\beta^\tau$ for the Curie–Weiss model. The classification of the points in $\tilde{\mathcal{E}}_{\beta,K}$ by type is done in Theorem 6.3 for $0 < \beta \leq \beta_c$ and $K > 0$, in which case $\tilde{\mathcal{E}}_{\beta,K}$ exhibits a continuous bifurcation, and in Theorem 6.4 for $\beta > \beta_c$ and $K > 0$, in which case $\tilde{\mathcal{E}}_{\beta,K}$ exhibits a discontinuous bifurcation. The associated limit theorems are given in Theorems 6.5 and 6.6. Except when $K = K_c^2(\beta)$ [Theorem 6.3(b)], the type of each of the global minimum points is 1. In these cases, the associated limit theorems are central-limit-type theorems with scalings $n^{1/2}$. When $K = K_c^2(\beta)$, we have $\tilde{\mathcal{E}}_{\beta,K} = \{0\}$, and the type of the minimum point at 0 is $r = 2$ or $r = 3$, depending on whether $0 < \beta < \beta_c$ or $\beta = \beta_c$. The associated limit theorems have noncentral-limit scalings $n^{3/4}$ or $n^{5/6}$, and in each case

$$P_{n,\beta,K}\{S_n/n^{1-1/2r} \in dx\} \quad \Longrightarrow \quad \text{const} \cdot \exp[-\text{const} \cdot x^{2r}]\, dx.$$

These nonclassical limit theorems signal the onset of a phase transition ([9], Section V.8). As $K$ increases through $K_c^2(\beta)$, the global minimum point at 0 bifurcates continuously into symmetric, nonzero global minimum points $\pm \tilde{z}(\beta, K)$.

We first consider $0 < \beta \leq \beta_c = \log 4$. According to Theorem 3.6, there exists a critical value

(6.5) $$K_c^{(2)}(\beta) = \frac{1}{2\beta c_\beta''(0)} = \frac{1}{4\beta e^{-\beta}} + \frac{1}{2\beta},$$

with the following properties:

- For $0 < K \leq K_c^{(2)}(\beta)$, $\tilde{\mathcal{E}}_{\beta,K} = \{0\}$.
- For $K > K_c^{(2)}(\beta)$, there exists $\tilde{z}(\beta, K) > 0$ such that $\tilde{\mathcal{E}}_{\beta,K} = \{\pm \tilde{z}(\beta, K)\}$.

The next theorem gives the type of each of these points in $\tilde{\mathcal{E}}_{\beta,K}$. The type is always 1 except when $K = K_c^{(2)}(\beta)$; in this case the global minimum point at 0 has type $r = 2$ if $0 < \beta < \beta_c$ and type $r = 3$ if $\beta = \beta_c$.

THEOREM 6.3. *Consider the BEG model, for which the canonical ensemble is given by* (3.2). *Let* $0 < \beta \leq \beta_c = \log 4$ *and define* $K_c^{(2)}(\beta)$ *by* (6.5). *The following conclusions hold:*

(a) *For* $0 < K < K_c^{(2)}(\beta)$, $z = 0$ *has type* $r = 1$.
(b) *Let* $K = K_c^{(2)}(\beta)$.

   (i) *For* $\beta < \beta_c$, $z = 0$ *has type* $r = 2$.
   (ii) *For* $\beta = \beta_c$, $z = 0$ *has type* $r = 3$.

(c) *For* $K > K_c^{(2)}(\beta)$ *and each choice of sign,* $z = \pm \tilde{z}(\beta, K)$ *has type* $r = 1$.

PROOF. (a) By (6.5), we have
$$G''_{\beta,K}(0) = 2\beta K(1 - 2\beta K c''_\beta(0))$$
$$= 2\beta K\left(1 - \frac{K}{K_c^{(2)}(\beta)}\right).$$

Therefore, $0 < K < K_c^{(2)}(\beta)$ implies that $G''_{\beta,K}(0) > 0$ and, thus, that $z = 0$ has type $r = 1$.

(b) For $K = K_c^{(2)}(\beta)$, $G''_{\beta,K}(0) = 0$. A simple calculation yields

(6.6)
$$G^{(4)}_{\beta,K}(0) = -(2\beta K)^4 c^{(4)}_\beta(0)$$
$$= -(2\beta K)^4 \cdot \frac{2e^{-\beta}(1 + 2e^{-\beta})(1 - 2e^{-\beta} - 8e^{-2\beta})}{(1 + e^{-\beta})^4}.$$

Therefore, for $\beta < \beta_c$, $G^{(4)}_{\beta,K}(0) > 0$ and for $\beta = \beta_c$, $G^{(4)}_{\beta,K}(0) = 0$. Computing the sixth derivative yields

(6.7) $$G^{(6)}_{\beta_c,K}(0) = 2 \cdot 3^4.$$

As a result, $z = 0$ has type 2 if $\beta < \beta_c$ and has type 3 if $\beta = \beta_c$.

(c) Lemma 3.7 states the existence and uniqueness of nonzero global minimum points $\pm \tilde{w}(\beta, K)$ of
$$F_{\beta,K}(w) = w^2/(4\beta K) - c_\beta(w) = G_{\beta,K}(w/(2\beta K)).$$

According to part (a) of the lemma, $F''_{\beta,K}(\tilde{w}(\beta, K)) > 0$. Since $\tilde{z}(\beta, K) = \tilde{w}(\beta, K)/(2\beta K)$, $F''_{\beta,K}(\tilde{w}(\beta, K)) > 0$ implies $G''_{\beta,K}(\tilde{z}(\beta, K)) > 0$. The symmetry of $G_{\beta,K}$ allows us to conclude that, for each choice of sign, $\pm \tilde{z}(\beta, K)$ has type $r = 1$. This completes the proof. □

We next classify by type the points in $\tilde{\mathcal{E}}_{\beta,K}$ for $\beta > \beta_c$ and $K > 0$. According to Theorem 3.8, there exists a critical value $K_c^{(1)}(\beta)$ with the following properties:

- For $0 < K < K_c^{(1)}(\beta)$, $\tilde{\mathcal{E}}_{\beta,K} = \{0\}$.
- For $K = K_c^{(1)}(\beta)$, there exists $\tilde{z}(\beta, K) > 0$ such that $\tilde{\mathcal{E}}_{\beta,K} = \{0, \pm \tilde{z}(\beta, K)\}$.
- For $K > K_c^{(1)}(\beta)$, $\tilde{\mathcal{E}}_{\beta,K} = \{\pm \tilde{z}(\beta, K)\}$.

The next theorem shows that the type of each of these points in $\mathcal{E}_{\beta,K}$ is 1.

THEOREM 6.4. *Consider the BEG model, for which the canonical ensemble is given by* (3.2). *Let $\beta > \beta_c$ and $K > 0$. The points in $\tilde{\mathcal{E}}_{\beta,K}$ all have type $r = 1$.*

PROOF. We first assume that $0 \in \tilde{\mathcal{E}}_{\beta,K}$, in which case $0 < K \leq K_c^{(1)}(\beta)$. Define $K_2 = 1/(2\beta c''_\beta(0))$. According to Theorem 3.8, we have $K_c^{(1)}(\beta) < K_2$. Since

$$
\begin{aligned}
G''_{\beta,K}(0) &= 2\beta K(1 - 2\beta K c''_\beta(0)) \\
&= 2\beta K \left(1 - \frac{K}{K_2}\right),
\end{aligned}
\tag{6.8}
$$

it follows that, whenever $0 < K \leq K_c^{(1)}(\beta)$, $1 > K/K_2$ and, thus, $G''_{\beta,K}(0) > 0$. We conclude that the global minimum point of $G_{\beta,K}$ at $z = 0$ has type $r = 1$, as claimed.

For $K \geq K_c^{(1)}(\beta)$, $\tilde{\mathcal{E}}_{\beta,K}$ also contains the symmetric, nonzero minimum points $\pm \tilde{z}(\beta, K)$ of $G_{\beta,K}$. Lemma 3.10 states the existence and uniqueness of nonzero global minimum points $\pm \tilde{w}(\beta, K)$ of

$$F_{\beta,K}(w) = w^2/(4\beta K) - c_\beta(w) = G_{\beta,K}(w/(2\beta K)).$$

Furthermore, according to part (a) of the lemma, $F''_{\beta,K}(\tilde{w}(\beta, K)) > 0$. Since $\tilde{z}(\beta, K) = \tilde{w}(\beta, K)/(2\beta K)$, $F''_{\beta,K}(\tilde{w}(\beta, K)) > 0$ implies $G''_{\beta,K}(\tilde{z}(\beta, K)) > 0$. The symmetry of $G_{\beta,K}$ allows us to conclude that, for each choice of sign, $\pm \tilde{z}(\beta, K)$ has type $r = 1$. This completes the proof. □

Theorems 6.1 and 6.2, together with the classification by type of the global minimum points of $G_{\beta,K}$, yield limit theorems for the $P_{n,\beta,K}$-distributions for appropriately scaled partial sums $S_n$ for the BEG model. The first, Theorem 6.5, states limits that are valid when $G_{\beta,K}$ has a unique global minimum point at $z = 0$. This is the case for $0 < \beta \leq \beta_c$, $0 < K \leq K_c^{(2)}(\beta)$ [Theorem 3.6(a)] and for $\beta > \beta_c$, $0 < K < K_c^{(1)}(\beta)$ [Theorem 3.8(a)]. The second, Theorem 6.6, states a law of large numbers and a conditioned limit that are valid when $G_{\beta,K}$ has multiple global minimum points.

In Theorem 6.5 $f_{0,\sigma^2(\beta,K)}$ denotes the density of a $N(0, \sigma^2(\beta, K))$ random variable with

$$\sigma^2(\beta, K) = \frac{2\beta K \cdot c''_\beta(0)}{G''_{\beta,K}(0)}. \tag{6.9}$$

When the type of the global minimum point at 0 is $r=1$, $\sigma^2(\beta, K) > 0$.

THEOREM 6.5. *Consider the BEG model, for which the canonical ensemble $P_{n,\beta,K}$ is given by (3.2). Suppose that $\tilde{\mathcal{E}}_{\beta,K} = \{0\}$ and let $r$ be the type of the point $z = 0$ as given in Theorems 6.3 and 6.4. The following conclusions hold:*

(a) $P_{n,\beta,K}\{S_n/n \in dx\} \Longrightarrow \delta_0$ *as* $n \to \infty$.
(b) *As* $n \to \infty$,

$$P_{n,\beta,K}\left\{\frac{S_n}{n^{1-1/2r}} \in dx\right\}$$

$$\Longrightarrow \begin{cases} f_{0,\sigma^2(\beta,K)}(x)\,dx, & \text{for } r=1, \\ \exp(-G^{(2r)}_{\beta,K}(0) \cdot x^{2r}/(2r)!)\,dx, & \text{for } r=2 \text{ or } r=3. \end{cases}$$

*When $r = 2$ [$K = K_c^{(2)}(\beta)$, $\beta < \beta_c$], $G^{(4)}_{\beta,K}(0)$ is given by (6.6), and when $r = 3$ [$K = K_c^{(2)}(\beta)$, $\beta = \beta_c$], $G^{(6)}_{\beta,K}(0) = 2 \cdot 3^4$.*

The last theorem states a law of large numbers and a conditioned limit that are valid when $G_{\beta,K}$ has multiple global minimum points. This holds in the following three cases:

1. $0 < \beta \leq \beta_c$ and $K > K_c^{(2)}(\beta)$, in which case the global minimum points are $\pm \tilde{z}(\beta, K)$ with $\tilde{z}(\beta, K) > 0$ [Theorem 3.6(b)];
2. $\beta > \beta_c$ and $K = K_c^{(1)}(\beta)$, in which case the global minimum points are $0$, $\pm \tilde{z}(\beta, K)$ with $\tilde{z}(\beta, K) > 0$ [Theorem 3.8(b)];
3. $\beta > \beta_c$ and $K > K_c^{(1)}(\beta)$, in which case the global minimum points are $\pm \tilde{z}(\beta, K)$ with $\tilde{z}(\beta, K) > 0$ [Theorem 3.8(c)].

In each case in which $G_{\beta,K}$ has multiple global minimum points, Theorems 6.3 and 6.4 states that all the global minimum points have type $r = 1$. For each global minimum point $\tilde{z}$ of type 1, we define the positive quantity

$$\sigma^2(\beta, K, z_j) = \frac{2\beta K \cdot c''_\beta(2\beta K z_j)}{G''_{\beta,K}(z_j)}. \tag{6.10}$$

THEOREM 6.6. *Consider the BEG model, for which the canonical ensemble $P_{n,\beta,K}$ is given by (3.2). Suppose that $\tilde{\mathcal{E}}_{\beta,K} = \{z_1, \ldots, z_m\}$ for $m = 2$ or $m = 3$. For each $j = 1, \ldots, m$, we define*

$$b_j = \frac{\sigma^2(\beta, z_j)}{\sum_{\ell=1}^m \sigma^2(\beta, z_\ell)},$$

where $\sigma^2(\beta, z_j)$ is the positive quantity defined in (6.10). Let $f_{0,\sigma^2(\beta,z_j)}$ be the density of a $N(0, \sigma^2(\beta, z_j))$ random variable. The following conclusions hold:

(a) $P_{n,\beta,K}\{S_n/n \in dx\} \Longrightarrow \sum_{j=1}^{m} b_j \delta_{z_j}$ as $n \to \infty$.
(b) There exists $\alpha = \alpha(z_j) > 0$ such that, for any $a \in (0, \alpha)$,

$$P_{n,\beta,K}\left\{\frac{S_n - nz_j}{n^{1/2}} \in dx \Big| \frac{S_n}{n} \in [z_j - a, z_j + a]\right\}$$
$$\Longrightarrow \quad f_{0,\sigma^2(\beta,z_j)}(x)\,dx \qquad as\ n \to \infty.$$

This completes our study of the limits for the $P_{n,\beta,K}$-distributions of appropriately scaled partial sums $S_n = \sum_{j=1}^{n} \omega_j$.

**Acknowledgments.** We thank Marius Costeniuc for supplying the proof of Proposition 3.4. We also thank one of the referees, on the basis of whose detailed comments we were able to improve this paper.

## REFERENCES

[1] AUSLOOS, M., CLIPPE, P., KOWALSKI, J. M. and PĘKALSKI, A. (1980). Magnetic lattice gas. *Phys. Rev. A* **22** 2218–2229.
[2] BARRÉ, J., MUKAMEL, D. and RUFFO, S. (2001). Inequivalence of ensembles in a system with long-range interactions. *Phys. Rev. Lett.* **87** 030601.
[3] BARRÉ, J., MUKAMEL, D. and RUFFO, S. (2002). Ensemble inequivalence in mean-field models of magnetism. *Dynamics and Thermodynamics of Systems with Long Range Interactions. Lecture Notes in Phys.* **602** 45–67. Springer, New York. MR2008178
[4] BLUME, M., EMERY, V. J. and GRIFFITHS, R. B. (1971). Ising model for the $\lambda$ transition and phase separation in $He^3$-$He^4$ mixtures. *Phys. Rev. A* **4** 1071–1077.
[5] CAPEL, H. W. (1966). On the possibility of first-order phase transitions in Ising systems of triplet ions with zero-field splitting. *Physica* **32** 966–987.
[6] CAPEL, H. W. (1967). On the possibility of first-order phase transitions in Ising systems of triplet ions with zero-field splitting II. *Physica* **33** 295–331.
[7] COSTENIUC, M., ELLIS, R. S. and TOUCHETTE, H. (2005). Complete analysis of phase transitions and ensemble equivalence for the Curie–Weiss–Potts model. *J. Math. Phys.* **46** 063301. MR2149837
[8] EISELE, T. and ELLIS, R. S. (1983). Symmetry breaking and random waves for magnetic systems on a circle. *Z. Wahrsch. Verw. Gebiete* **63** 297–348. MR705628
[9] ELLIS, R. S. (1985). *Entropy. Large Deviations and Statistical Mechanics*. Springer, New York. MR793553
[10] ELLIS, R. S., HAVEN, K. and TURKINGTON, B. (2000). Large deviation principles and complete equivalence and nonequivalence results for pure and mixed ensembles. *J. Statist. Phys.* **101** 999–1064. MR1806714
[11] ELLIS, R. S., MONROE, J. L. and NEWMAN, C. M. (1976). The GHS and other correlation inequalities for a class of even ferromagnets. *Comm. Math. Phys.* **46** 167–182. MR395659

[12] Ellis, R. S. and Newman, C. M. (1979). Limit theorems for sums of dependent random variables occurring in statistical mechanics. *Z. Wahrsch. Verw. Gebiete* **44** 117–139. MR503333

[13] Ellis, R. S., Newman, C. M. and O'Connell, M. R. (1981). The GHS inequality for a large external field. *J. Statist. Phys.* **26** 37–52. MR643702

[14] Ellis, R. S., Newman, C. M. and Rosen, J. S. (1980). Limit theorems for sums of dependent random variables occurring in statistical mechanics, II: Conditioning, multiple phases, and metastability. *Z. Wahrsch. Verw. Gebiete* **51** 153–169. MR566313

[15] Ellis, R. S., Touchette, H. and Turkington, B. (2004). Thermodynamic verses statistical nonequivalence of ensembles for the mean-field Blume–Emery–Griffiths model. *Phys. A* **335** 518–538. MR2044158

[16] Hoston, W. and Berker, A. N. (1991). Multicritical phase diagrams of the Blume–Emery–Griffiths model with repulsive biquadratic coupling. *Phys. Rev. Lett.* **67** 1027–1030.

[17] Kivelson, S. A., Emery, V. J. and Lin, H. Q. (1990). Doped antiferromagnets in the weak-hopping limit. *Phys. Rev. B* **42** 6523–6530.

[18] Lajzerowicz, J. and Sivardière, J. (1975). Spin-1 lattice-gas model. I. Condensation and solidification of a simple fluid. *Phys. Rev. A* **11** 2079–2089.

[19] Newman, K. E. and Dow, J. D. (1983). Zinc-blende-diamond order-disorder transition in metastable crystalline $(GzAs)_{1-x}Ge_{2x}$ alloys. *Phys. Rev. B* **27** 7495–7508.

[20] Nienhuis, B., Berker, A. N., Riedel, E. K. and Schick, M. (1979). First- and second-order phase transitions in Potts models: Renormalization-group solution. *Phys. Rev. Lett.* **43** 737–740.

[21] Otto, P. T. (2004). Study of equilibrium macrostates for two models in statistical mechanics. Ph.D. dissertation, Univ. Massachusetts, Amherst. Available at http://www.math.umass. edu/~rsellis/pdf-files/Otto-thesis.pdf.

[22] Rockafellar, R. T. (1970). *Convex Analysis*. Princeton Univ. Press. MR274683

[23] Schick, M. and Shih, W.-H. (1986). Spin-1 model of a microemulsion. *Phys. Rev. B* **34** 1797–1801.

[24] Sivardière, J. and Lajzerowicz, J. (1975). Spin-1 lattice-gas model. II. Condensation and phase separation in a binary fluid. *Phys. Rev. A* **11** 2090–2100.

[25] Sivardière, J. and Lajzerowicz, J. (1975). Spin-1 lattice-gas model. III. Tricritical points in binary and ternary fluids. *Phys. Rev. A* **11** 2101–2110.

R. S. Ellis  
Department of Mathematics and Statistics  
University of Massachusetts  
Amherst, Massachusetts 01003  
USA  
e-mail: rsellis@math.umass.edu

P. T. Otto  
Department of Mathematics  
Union College  
Schenectady, New York 12308  
USA  
e-mail: ottop@union.edu

H. Touchette  
School of Mathematical Sciences  
Queen Mary, University of London  
London E1 4NS  
United Kingdom  
e-mail: htouchet@alum.mit.edu